# ON THE IRRELEVANT DISORDER REGIME OF PINNING MODELS[1]

By Giambattista Giacomin and Fabio Lucio Toninelli

*Université Paris Diderot and ENS Lyon*

Recent results have lead to substantial progress in understanding the role of disorder in the (de)localization transition of polymer pinning models. Notably, there is an understanding of the crucial issue of *disorder relevance* and *irrelevance* that is now rigorous. In this work, we exploit interpolation and replica coupling methods to obtain sharper results on the irrelevant disorder regime of pinning models. In particular, in this regime, we compute the first order term in the expansion of the free energy close to criticality and this term coincides with the first order of the formal expansion obtained by field theory methods. We also show that the quenched and quenched averaged correlation length exponents coincide, while, in general, they are expected to be different. Interpolation and replica coupling methods in this class of models naturally lead to studying the behavior of the intersection of certain renewal sequences and one of the main tools in this work is precisely renewal theory and the study of these intersection renewals.

**1. Introduction.** The role played by quenched disorder in statistical mechanics models is still little understood, not only from the mathematical standpoint, but also at less rigorous levels of analysis. Still, some physical approaches, in spite of being nonrigorous, provide predictions that are very intriguing for mathematicians, for at least two reasons. First, disordered systems are *doubly probabilistic*, with two sources of randomness that enter at the same time, but in very distinct ways, and this interplay has clearly demanded wholly new ideas that have then played a role well beyond the realm of statistical mechanics (we cite here the particularly remarkable example of all the mathematical tools that have been developed around spin glasses [7, 8, 25, 28]). Second, *solutions* or *conjectures*, set forth mostly in

Received September 2007.
[1]Supported in part by ANR project POLINTBIO.

theoretical physics, are often very bold and have a very novel character, so classical approaches are insufficient to prove (or disprove) them. Also, for this aspect, spin glasses possibly constitute the main example, but one should not neglect the many other extremely relevant models, some of these already developed in the references we have mentioned (e.g., random walks and diffusions in random environments).

Polymer models with random potentials fully fall into this category and the list of references would be practically endless, partly because of the natural link between polymer models and other models involving, in particular, phase boundaries and interfaces. What we are going to consider, the pinning models, constitute a subclass of polymers with random potentials. Pinning models are directed polymer models which are particularly attractive, both for their very large spectrum of applications and for the fact that they are fully *exactly solvable* in their nondisordered (also called *homogeneous* or *pure*) version, yet exhibit, for example, phase transitions of all possible orders [15]. Moreover, from a more strictly mathematical viewpoint, the close link between pinning models and renewal theory is another matter of interest since nontrivial *cross-questions* arise (we cite, e.g., [10] and [19], but the present work itself is heavily based on the interplay with renewal theory).

Pinning models are indexed by a positive parameter that we call $\alpha$ (in agreement with the probabilistic literature on heavy-tailed distributions and somewhat in disagreement with the physical literature where $\alpha$ is one of the critical exponents). Varying $\alpha$, one *walks through* a number of different models (pinning of polymers on defect lines in different dimensions, the Poland–Scheraga model of DNA denaturation, Wetting models, etc. [15, 18]) and, at the same time, a variety of critical behaviors emerge, notably, as pointed out above, phase transitions of any order, at least if one looks only at nondisordered models. A very intriguing question therefore arises: what is the effect of disorder (even a tiny amount of it) on the transition?

This question is highly debated and certainly not just for pinning models. Very well known is the result of Aizenman and Wehr [1] on the regularizing effect of disorder on the transition of certain models (in particular, Ising models) that makes rigorous an argument due to Imry and Ma [24] which is based on the comparison between the competing effects of boundary conditions and of random field fluctuations. Their method extends to $(1 + 2)$-dimensional interfaces [9], but not to pinning models. A variety of authors (e.g., [12, 16]) have applied to pinning models a nonrigorous argument originally due to A. B. Harris predicting that introducing a small disorder leads to a change of the critical properties of the model (normally, in this case, one says that the disorder is *relevant*) when the *specific heat* diverges at criticality in the corresponding pure model. For pinning models, the Harris criterion boils down to predicting that disorder is relevant for $\alpha > 1/2$ and irrelevant if $\alpha < 1/2$, with substantial uncertainties in the physical literature for the marginal case $\alpha = 1/2$ (see, however, [20]). Moreover, for pinning models, the theoretical physics arguments suggest that the critical point of the quenched model coincides, when the disorder is irrelevant, with the critical point of the annealed model, which is nothing but a homogeneous model.

Thanks to the work of Alexander [3] for the irrelevant disorder regime and to the work of the authors [22] when disorder is relevant, we now know that the Harris criterion yields correct predictions. In spite of the fact that several open questions still stand, this is quite a satisfactory situation that provides motivation to investigate further. In particular, one of the two authors has recently found a different approach [30] to Alexander's result, based on the *replica coupling method* developed in [23], in the context of spin glasses.

In this work, we are going to exploit the interpolation method further, making rigorous some physical predictions when disorder is irrelevant. For instance, we identify the first order term in the expansion of the free energy at small coupling close to criticality; see Theorem 2.3 below. Moreover, we establish some results that, while possibly intuitive in the rationale of the irrelevant disorder regime, were not clear cut in the physics literature. We study, for example, the correlation length exponent in the quenched and quenched averaged setup. In general, these two exponents are not expected to coincide (we cite, e.g., the quantum Ising chain with random transverse field studied in [14]). In this work, we show that the two exponents coincide when disorder is irrelevant.

Another important aspect of the approach we take is that by studying *replicas* of the original systems, challenging questions on the intersection of independent renewal sequences come up. This is a particular instance of the much studied problem of intersection of independent Markov processes (recall that renewal processes are just random walks with positive increments). It is certainly not possible to give a proper account of this activity and we limit ourselves to mentioning [5], and references therein, in which is pointed out, in particular, the difficulty of the task of determining the law of intersection renewals. In our case, the renewal sequences involved have exponential inter-arrival tails and this presents challenges that are very different from the ones encountered when dealing with the more widely studied case of heavy-tail inter-arrivals laws [5]. In this work, we tackle only the questions that are relevant to our analysis, but we feel that some of the results may be of independent interest and that replica and interpolation methods may lead to further novel and nontrivial questions on renewal sequences.

## 2. The general (disordered) model and its homogeneous counterpart.

2.1. *The basic renewal sequence.* Let $\tau := \{\tau_n\}_{n=0,1,\dots}$ be a homogeneous (nondelayed) [4] renewal sequence, namely $\tau$ is a sequence of random variables with $\tau_0 = 0$ and, for $n \geq 1$, $\tau_n = \sum_{j=1}^n \eta_j$, where $\{\eta_n\}_{n=1,2,\dots}$ is an

i.i.d. sequence. Note that the law of $\eta_n$ coincides with the law of $\tau_1$. We assume that $\tau_1$ takes values in $\mathbb{N} \cup \{\infty\} = \{1, 2, \dots\} \cup \{\infty\}$ and we set $K(n) := \mathbf{P}(\tau_1 = n)$. We call $K(\cdot)$ an *inter-arrival law* and, in general, we are going to refer to a renewal sequence with inter-arrival law $F(\cdot)$ as $F(\cdot)$-renewal. Throughout the paper, we assume that

$$(2.1) \qquad\qquad K(n) = \frac{L(n)}{n^{1+\alpha}},$$

where $\alpha > 0$ and $L(\cdot)$ is a slowly varying function. We recall that $L : (0, \infty) \to (0, \infty)$ is slowly varying if it is measurable and if $\lim_{x\to\infty} L(cx)/L(x) = 1$ for every $c > 0$, where the notation $a(x) \overset{x\to l}{\sim} b(x)$ means that $\lim_{x\to l} a(x)/b(x) = 1$. We refer to [6] for the full theory of *slow variation*. For the sake of generality, we are not going to assume that $\sum_{n\in\mathbb{N}} K(n) = 1$, but rather that $\sum_{n\in\mathbb{N}} K(n) \leq 1$, and we set $K(\infty) := 1 - \sum_{n\in\mathbb{N}} K(n)$ (see, however, Remark 2.1 below). This means that, in general, $\tau_1$ takes on finite values only with probability $1 - K(\infty)$. When $K(\infty) > 0$, we say that the renewal is terminating or transient, while if $K(\infty) = 0$, we say that it is persistent or recurrent.

For notational convenience, we are also going to look at $\tau$ as a subset of $\mathbb{N} \cup \{0, \infty\}$. Note that in the terminating case, almost surely $\tau$ is a finite set (containing $\infty$) and in the persistent case, $\tau$ contains infinitely many points (but $\infty \notin \tau$). As is customary, the function $n \mapsto \mathbf{P}(n \in \tau)$, $\tau$ a general renewal sequence, is called the *renewal function* of $\tau$. The renewal function is related to the inter-arrival law by the recurrence scheme

$$(2.2) \qquad \mathbf{P}(n \in \tau) = \mathbf{1}_{\{0\}}(n) + \sum_{j=1}^{n} \mathbf{P}(\tau_1 = j)\mathbf{P}(n - j \in \tau).$$

We cite here the following important consequence of the sharp renewal estimates proven in [13, 17]: if $\tau$ is a persistent $K(\cdot)$-renewal, then, for every $\alpha \in (0, 1)$,

$$(2.3) \qquad \mathbf{P}(n \in \tau) \overset{n\to\infty}{\sim} \frac{\alpha \sin(\pi\alpha)}{\pi} \frac{1}{L(n)n^{1-\alpha}}.$$

**2.2. *The general model.*** The disordered pinning model of finite size $N \in \mathbb{N}$ and parameters $\beta \geq 0$ and $h \in \mathbb{R}$ is defined by introducing the new probability measure $\mathbf{P}_{N,\omega}(= \mathbf{P}_{N,\omega,\beta,h})$ and the (realization of the) sequence of independent standard normal random variables $\omega := \{\omega_n\}_{n\in\mathbb{N}}$ via the formula

$$(2.4) \qquad \frac{d\mathbf{P}_{N,\omega}}{d\mathbf{P}}(\tau) := \frac{1}{Z_{N,\omega}} \exp\left(\sum_{n=1}^{N} (\beta\omega_n + h)\mathbf{1}_{n\in\tau}\right)\mathbf{1}_{N\in\tau},$$

where $Z_{N,\omega}$ is the normalization (or *partition*) function. All major results on this model are stated assuming that $K(\infty) = 0$, without loss of generality (see Remark 2.1 below), but terminating renewals play a central role in the technical arguments.

Informally, $\mathbf{P}_{N,\omega}$ describes a point process that favors trajectories, that is, random subsets of the integer numbers, which maximize the *energy*, that is, the sum over the subset of the (inhomogeneous) quantity $\beta\omega_n + h$. Due to the fact that values of the energy close to the maximal one are typically reached only by few configurations, a nontrivial energy-entropy competition arises. This model presents a localization/delocalization transition, in the sense that if *overall* the energy contributions are negative, for example, if $h$ is negative and large in absolute value, then, in the limit as $N \to \infty$, the process trajectories concentrate on sets containing only a few points and these points are close to the boundary of the system (this is what we call a *delocalization* phenomenon). The complementary situation is observed when the energy contributions are *overall* positive and the size of the random sets that are typically observed for $N$ large is a positive fraction of $N$ (*localization*). We are being very imprecise about what we mean by *overall positive* or *overall negative*; in a sense, making this concept precise is the central issue in comprehending the (de)localization transition. We refer, for example, to [18], Chapter 1, for a substantially more detailed discussion of the model, an overview of the literature and a survey of the very many contexts in which this model has been proposed and studied.

The law of $\omega$ is denoted by $\mathbb{P}$. The model has to be understood as a *quenched* model, but, in this work, we will also focus on the *quenched averaged* measure $\mathbb{E}\mathbf{P}_{N,\omega}$ (not to be confused with the *annealed* measure that we discuss below). Quenched and quenched averaged quantities may coincide in the limit as $N \to \infty$ (this is the phenomenon known as *self-averaging*). The first and (possibly) most important of the self-averaging quantities is the free energy

$$(2.5) \qquad \mathrm{F}(\beta, h) := \lim_{N \to \infty} \frac{1}{N} \log Z_{N,\omega}.$$

The limit is to be understood in the $\mathbb{P}(\mathrm{d}\omega)$-almost sure or in the $L^1(\mathbb{P}(\mathrm{d}\omega))$ sense. The existence of such a limit follows by standard arguments (see, e.g. [18], Chapter 4). It is also standard to show that $\mathrm{F}(\beta, h) \geq 0$ and to split the parameter space into a *delocalized* region $\mathcal{D} := \{(\beta, h) : \mathrm{F}(\beta, h) = 0\}$ and a *localized* one $\mathcal{L} := \{(\beta, h) : \mathrm{F}(\beta, h) > 0\} = \{(\beta, h) : h > h_c(\beta)\}$ with $h_c(\beta) := \inf\{h : \mathrm{F}(\beta, h) > 0\}$ [18], Chapter 1.

We quickly recall that this splitting of the phase space into localized and delocalized regions corresponds to sharply different path properties in the limit of large values of $N$ (see [18], Chapter 7, and references therein). In particular, in [21], it is shown that the weak limit $\mathbf{P}_{\infty,\omega}$ of $\{\mathbf{P}_{N,\omega}\}_N$ exists $\mathbb{P}(\mathrm{d}\omega)$-a.s. for $h > h_c(\beta)$. Moreover, in [21], it is shown that $\mathtt{F}(\cdot, \cdot)$ is $C^\infty$ for $h > h_c(\beta)$. We will also occasionally need the definition

$$(2.6) \qquad \mathtt{F}_N(\beta, h) := \frac{1}{N}\mathbb{E}\log Z_{N\omega}(\beta, h).$$

2.3. *The homogeneous model.* With abuse of notation, we let $\mathbf{P}_{N,h} = \mathbf{P}_{N,\omega,0,h}$. This case is, of course, much easier to handle since inhomogeneous potentials are no longer present. As a matter of fact, the model is completely solvable. Understanding the homogeneous model in some detail is quite central for this work and we therefore give a detailed sketch of what is known about the model. More details can be found in [18], Chapters 1 and 2.

With further abuse of notation, given $h$, we denote by $\mathtt{F}[=\mathtt{F}(h)]$ the unique solution to

$$(2.7) \qquad \sum_{n\in\mathbb{N}} K(n)\exp(-\mathtt{F}n) = \exp(-h),$$

if such a solution exists, and $\mathtt{F} := 0$ otherwise. Note that if $\mathtt{F} < 0$, by (2.1), the left-hand side of (2.7) is equal to $\infty$ and from this, one easily infers that a solution $\mathtt{F}$ to (2.7) can be found only if $\sum_{n\in\mathbb{N}} K(n) \le \exp(-h)$, that is, if $h \ge h_c := -\log\sum_{n\in\mathbb{N}} K(n)$. Elaborating slightly further, one also easily sees that $h \mapsto \mathtt{F}(h)$ is nondecreasing and is increasing if $h > h_c$. Moreover, it is a continuous function; in fact, it is real analytic, except at $h = h_c$. It is not difficult to show that the following identity holds for every $N$:

$$(2.8) \qquad Z_{N,h} = \exp(\mathtt{F}(h)N)\mathbf{P}_h(N \in \tau),$$

where, under $\mathbf{P}_h$, the sequence $\tau$ is a renewal sequence with inter-arrival law $n \mapsto \mathbf{P}_h(\tau_1 = n) = \exp(h)K(n)\exp(-\mathtt{F}(h)n)$. By the very definition of $\mathtt{F}(h)$, such a renewal sequence is terminating if $h < h_c$ and is persistent if $h \ge h_c$. From (2.8), it directly follows that

$$(2.9) \qquad \mathtt{F}(h) = \lim_{N\to\infty}\frac{1}{N}\log Z_{N,h},$$

so $\mathtt{F}(h)$ is the free energy of the model and, going back to (2.5), we note that $\mathtt{F}(h) = \mathtt{F}(0, h)$ and $h_c = h_c(0)$.

REMARK 2.1. The explicit solution we have outlined shows, in particular, that there is no loss of generality in assuming that $\sum_{n\in\mathbb{N}} K(n) = 1$ since this simply leads to a change in the value of $h_c$. The same is also true for the disordered case introduced in Section 2.2 (for a detailed discussion of this point, we refer to [18], Chapter 1). We are therefore going to make this assumption throughout the remainder of the paper, so $h_c = h_c(0) = 0$.

REMARK 2.2. Note that $\mathtt{F}(h) \geq 0$ is a consequence of the subexponential character of $K(\cdot)$: from (2.7), one sees that $\mathtt{F}(\cdot) \geq 0$ if and only if $\sum_{n \in \mathbb{N}} K(n) \times \exp(\delta n) = \infty$ for every $\delta > 0$. If this latter condition fails, then there exists some $h$ such that $\mathtt{F}(h) < 0$, but (2.8) still holds. The technical arguments in this work at times rely on pinning problems for homogeneous exponentially decaying inter-arrival laws [i.e., such that $\sum_{n \in \mathbb{N}} K(n) \exp(\delta n) < \infty$ for some $\delta > 0$] and, even if we are not interested in the negative free energy regime, we wish to point out that the solution scheme for these models is the same as that for the subexponentially decaying inter-arrivals.

2.4. *Quenched and annealed models.* Important bounds for quenched systems come from the *annealing* procedure. This simply involves exchanging the order of the disorder average and the logarithm in the definition of the quenched averaged model, namely, at the partition function level: $\mathbb{E} \log Z_{N,\omega} \leq \log \mathbb{E} Z_{N,\omega}$. The latter expression is nothing but the logarithm of the partition function of an homogeneous model with potential equal to $h + (\beta^2/2)$. Such a model is referred to as *annealed* and its free energy is therefore simply $\mathtt{F}(0, h + (\beta^2/2))$, which coincides with $\mathtt{F}(h + (\beta^2/2))$ in the shorthand notation of Section 2.3. Of course, the critical point for the annealed model is just $h_c^a(\beta) := -(\beta^2/2)$ [here, we use the fact that $h_c = h_c(0) = 0$ with the assumption that $K(\infty) = 0$]. Therefore, for every $\Delta \in \mathbb{R}$, we have

$$\mathtt{F}(\beta, h_c^a(\beta) + \Delta) \leq \mathtt{F}(0, \Delta) \tag{2.10}$$

and $h_c(\beta) \geq h_c^a(\beta)$ for every $\beta$.

Quenched-to-annealed comparisons are generally strict, but it has been recently proven (in [3]) that if either

$$\alpha \in (0, 1/2) \quad \text{or} \quad \left\{ \alpha = 1/2 \text{ and } \sum_n \frac{1}{nL(n)^2} < \infty \right\}, \tag{2.11}$$

then there exists some $\beta_0 > 0$ such that

$$h_c(\beta) = h_c^a(\beta) \qquad \text{if } \beta \leq \beta_0, \tag{2.12}$$

as a consequence of the fact that for every $\varepsilon > 0$, there exists some $\Delta_0 > 0$ such that

$$\mathtt{F}(\beta, h_c^a(\beta) + \Delta) \geq (1 - \varepsilon)\mathtt{F}(0, \Delta) \tag{2.13}$$

for $\Delta \leq \Delta_0$ and $\beta \leq \beta_0$. The proof in [3] is based on a modified second moment method. The alternative proof in [30] is instead based on the interpolation method and it is precisely by exploiting this second method further that we are able to sharpen (2.13). In order to better understand the results that we have just mentioned, as well as the role of (2.11), we now make a brief detour that is also going to allow us to state a first new result.

2.5. *Harris criterion and renewal intersections.* A key question in several instances is understanding whether or not introducing the disorder *substantially* changes the behavior of a model. This question is particularly intriguing close to criticality and, as mentioned in Section 1, a somewhat general argument to determine whether the disorder is relevant or irrelevant has been proposed by A. B. Harris. In particular, arguing à la Harris for disordered pinning models, one obtains that for $\alpha < 1/2$, quenched and annealed criticality points, as well as quenched and annealed critical properties, should coincide for small values of $\beta$. Note that (2.12) and (2.13) prove the correctness of the physical claims for $\alpha < 1/2$, that is, when disorder is irrelevant. Actually, the small disorder expansion arguments in [16] suggest that when disorder is irrelevant, one should be able to write

$$(2.14) \qquad \mathrm{F}(\beta, h_c^a(\beta) + \Delta) = \mathrm{F}(0, \Delta) - \frac{\beta^2}{2}(\partial_\Delta \mathrm{F}(0, \Delta))^2 + \cdots$$

for $\beta$ small and $\Delta \searrow 0$. For such an expansion, it is of course crucial that $(\partial_\Delta \mathrm{F}(0, \Delta))^2 = o(\mathrm{F}(0, \Delta))$ as $\Delta \searrow 0$ and this will, in fact, be shown to hold under assumption (2.11). In fact, (2.14) itself can actually be made rigorous and an upper bound corresponding to (2.14) has already been proven in [30]. Here, we are going to prove the opposite bound. More precisely, we establish the following:

THEOREM 2.3. *If (2.11) holds, then*

$$(2.15) \qquad \lim_{\beta \searrow 0} \limsup_{\Delta \searrow 0} \left| \frac{\mathrm{F}(0, \Delta) - \mathrm{F}(\beta, h_c^a(\beta) + \Delta)}{(\beta^2/2)(\partial_\Delta \mathrm{F}(0, \Delta))^2} - 1 \right| = 0.$$

What we are going to prove goes beyond the content of Theorem 2.3, in the sense that estimates for $\beta$ and $\Delta$ small but finite are also established. In particular, it is natural to ask for which values of $\beta$ the small-$\Delta$ expansion in (2.14) can be performed. Alternatively, one can ask what the value of $\beta_0$ is in (2.12). This is what we are going to explain next and this is also going to clarify the role of the hypothesis (2.11) in probabilistic terms.

A crucial mathematical object that comes up in [3, 30] and (more implicitly) in [12, 16] is the *intersection renewal* $\tau \cap \tau'$, where $\tau'$ is an independent copy of $\tau$. It is well known and straightforward to show that if $\tau$ and $\tau'$ are two independent (general) renewals, then the sequence (or random set) $\tau \cap \tau'$ is also a renewal and $\mathbf{P}(n \in \tau \cap \tau') = \mathbf{P}(n \in \tau)\mathbf{P}(n \in \tau')$, so the renewal function is explicit and, by using (2.2), one can then extract the inter-arrival law of $\tau \cap \tau'$. We point out that the implicit character of (2.2) means that this procedure is of nonimmediate applicability. There are, however, some properties of $\tau \cap \tau'$ that one can easily address, in particular, whether $\tau \cap \tau'$ is terminating or persistent. Restricting to the case in which both $\tau$ and $\tau'$ are $K(\cdot)$-renewals, by using the identity $\mathbf{E}[|\tau \cap \tau'|] = \sum_n \mathbf{P}(n \in \tau)^2$ and (2.3), one sees that $\tau \cap \tau'$ is almost surely a finite set and therefore the renewal is terminating (in spite of the fact that $\tau$ and $\tau'$ are persistent) if (2.11) holds. Actually, since, in full generality, the expected size of a renewal set coincides with the reciprocal of the probability that the inter-arrival variable takes the value $+\infty$ (e.g., [6], Section 8.7), (2.11) is necessary and sufficient for the intersection of two independent $K(\cdot)$-renewals to be terminating.

If one now considers the 2-replica homogeneous pinning model with free energy equal to

$$(2.16) \qquad \lim_{N \to \infty} \frac{1}{N} \log \mathbf{E}^{\otimes 2}\left[ \exp\left( \lambda \sum_{n=1}^{N} \mathbf{1}_{n \in \tau \cap \tau'} \right) ; N \in \tau \cap \tau' \right],$$

the computation of such a limit falls in the realm of the general theory recalled in Section 2.3. In particular, since $\tau \cap \tau'$ is terminating and the inter-arrival law of $\tau \cap \tau'$ is subexponential [in fact, $\mathbf{P}^{\otimes 2}(\inf\{n > 0 : n \in \tau \cap \tau'\} = N) \geq K(N)^2$], by applying (2.7)–(2.9), one sees that the expression in (2.16) is zero for $\lambda \leq \lambda_0$, with $\lambda_0 := -\log(1 - \mathbf{P}^{\otimes 2}((\tau \cap \tau') = \{0\})) > 0$. The quantity $\beta_0$ mentioned above [see (2.12)] may be chosen to be equal to $\lambda_0/2$.

Besides the quantitative estimate, the aim of what we have just explained is to emphasize that the irrelevant disorder regime holds when the renewal intersection built from the $K(\cdot)$-renewals is terminating (delocalized) and the coupling parameter is so small that it does not make them persistent (localized).

2.6. *On the irrelevant disorder regime.* We are now going to state more results in the spirit of Theorem 2.3. In order to do this, we need to recall a quantity introduced in [2, 21]:

$$(2.17) \qquad \mu_N(\beta, h) := -\frac{1}{N} \log \mathbb{E}\left[ \frac{1}{Z_{N,\omega}} \right].$$

We let $\mu(\beta, h) := \lim_{N \to \infty} \mu_N(\beta, h)$ (the limit exists by superadditivity).

A certain number of facts are known about $\mu$. First of all, we have, in general, that $0 \leq \mu \leq \mathtt{F}$, the lower bound being a consequence of the subexponential character of $K(\cdot)$ and the upper one of Jensen's inequality. However, a stronger statement was proven in [21]: if $h > h_c(\beta)$, then

$$(2.18) \qquad 0 < \mu(\beta, h) < \mathtt{F}(\beta, h).$$

In particular, $\mu(\beta, h) > 0$ if and only if $(\beta, h) \in \mathcal{L}$. We will see in a moment that a question of great interest is whether or not $\mu$ and $\mathtt{F}$ have the same critical behavior close to $h_c(\beta)$. For the moment, let us mention that in [31], it was proven that

$$(2.19) \qquad c(\beta) \frac{\mathtt{F}(\beta, h)^2}{\partial_h \mathtt{F}(\beta, h)} < \mu(\beta, h)$$

for some positive $c(\beta)$ if, say, $0 < h - h_c(\beta) \le 1$. This can be expressed in a less precise, but more intuitive, way: assume that for $h \searrow h_c(\beta)$,

$$
\begin{aligned}
(2.20) \qquad & \textsc{f}(\beta, h) \sim c_\textsc{f}(\beta)(h - h_c(\beta))^{\nu_\textsc{f}}, \\
& \mu(\beta, h) \sim c_\mu(\beta)(h - h_c(\beta))^{\nu_\mu}.
\end{aligned}
$$

Then, (2.19) and the upper bound in (2.18) imply that $\nu_\textsc{f} \le \nu_\mu \le \nu_\textsc{f} + 1$. Our Theorem 2.4 below says, in particular, that in the irrelevant disorder regime, $\nu_\textsc{f} = \nu_\mu$.

Before we formulate Theorem 2.4, it is useful to discuss why $\mu$ is an interesting quantity to look at and why we wish to compare the critical behavior of $\textsc{f}$ and $\mu$.

We demonstrate the relevance of $\mu$ by giving three instances in which it appears:

1. Let $\Delta_N$ be the largest gap in the renewal up to $N$, that is,

$$
(2.21) \qquad \Delta_N := \max\{j - i : 1 \le i < j \le N \text{ and } \tau \cap \{i, \ldots, j\} = \varnothing\}.
$$

Then, for $h > h_c(\beta)$, one has, for every $\varepsilon > 0$,

$$
(2.22) \qquad \lim_{N \to \infty} \mathbf{P}_{N,\omega}\left( \frac{1 - \varepsilon}{\mu(\beta, h)} \le \frac{\Delta_N}{\log N} \le \frac{1 + \varepsilon}{\mu(\beta, h)} \right) = 1
$$

in $\mathbb{P}(\mathrm{d}\omega)$-probability.

This was proven in [2] in a related model (the *copolymer at a selective interface*) and, in the present context, in [21], Theorem 2.5.

2. Let $\Delta_x := \tau_{\iota(x)+1} - \tau_{\iota(x)}$ with $\iota(x)$ equal to the value of $j$ for which $x \in \{\tau_j, \tau_j + 1, \ldots, \tau_{j+1} - 1\}$. For every $\varepsilon > 0$, there exists $c = c(\varepsilon, \beta, h) > 0$ such that for every $x < N$,

$$
(2.23) \qquad c^{-1} e^{-\mu(\beta,h)(1+\varepsilon)n} \le \mathbb{E}\mathbf{P}_{N,\omega}(\Delta_x = n) \le c e^{-\mu(\beta,h)(1-\varepsilon)n},
$$

where the lower bound holds for $n \le N/2$ and the upper one for $n \in \mathbb{N}$. This can be proven in analogy with [21], Proposition 2.4 and it is detailed in [18], Chapter 7.

3. For the system in the localized phase, two distinct *correlation lengths* are naturally defined. One starts from the two-point function defined as

$$
(2.24) \qquad C_\omega(k) := \mathbf{P}_{\infty,\omega}(0 \in \tau, k \in \tau) - \mathbf{P}_{\infty,\omega}(0 \in \tau)\mathbf{P}_{\infty,\omega}(k \in \tau),
$$

where $\mathbf{P}_{\infty,\omega}$ is the bi-infinite volume measure built in a natural way in the localized regime, with $\omega$ an i.i.d. bi-infinite sequence (see [21] for details). In [21], it is shown that $_\omega(k)$ vanishes exponentially for $k \to \infty$ and the same holds [$\mathbb{P}(\mathrm{d}\omega)$-almost surely] without taking the disorder average. It is then natural to call the inverse of the rate of exponential decay of $\mathbb{E}C_\omega(k)$ [resp. of $C_\omega(k)$] the *average correlation length* $\xi_{\mathrm{av}}$ (resp., *typical correlation length* $\xi_{\mathrm{typ}}$).

From general principles, one expects that close to criticality, $\xi_{\mathrm{av}}$ and $\xi_{\mathrm{typ}}$ diverge like $1/\mu$ and $1/\mathsf{F}$, respectively, and in some specific examples, this can be proven. This is, in particular, true when $K(n) = \mathbf{P}(\inf\{k : S_k = 0\} = 2n)$ and $\{S_n\}_{n \geq 0}$ is the one-dimensional simple random walk started from 0. In this case, it was proven in [31] that for every $h > h_c(\beta)$,

$$(2.25) \qquad \xi_{\mathrm{av}}(\beta, h) = \mu(\beta, h)^{-1} \qquad \text{while } \xi_{\mathrm{typ}}(\beta, h) = \mathsf{F}(\beta, h)^{-1}.$$

Note that, for this choice of $K(\cdot)$, one has $\alpha = 1/2$ and $L(\cdot)$ asymptotically constant in (2.1), so (2.11) is *not* satisfied.

For a general inter-arrival law $K(\cdot)$ satisfying (2.1), (2.25) does not hold since this is already the case in the nondisordered setup [19]. However, in the nondisordered setup, (2.25) does hold for $h - h_c$ sufficiently small and a (possibly weaker) form of (2.25) is expected for the whole class of disordered models sufficiently close to criticality (see [32] for some results in this direction).

In general, there is no reason to believe that $\xi_{\mathrm{typ}}$ and $\xi_{\mathrm{av}}$ have the same critical behavior (see, e.g., [14]).

As we have already mentioned, in the irrelevant disorder regime, we can prove that $\mu$ and $\mathsf{F}$ behave in essentially the same way close to criticality (in particular, the critical exponents coincide) and for both, we have a control over the first-order term in the small-disorder expansion near criticality [cf. (2.14)].

THEOREM 2.4. *If (2.11) holds, then there exist positive constants $\beta_0$, $C$ and $\Delta_0$ such that, for $\beta \leq \beta_0$ and $0 \leq \Delta \leq \Delta_0$,*

$$(2.26) \quad \mathsf{F}(0, \Delta) - 9\beta^2(\partial_\Delta \mathsf{F}(0, \Delta))^2 \leq \mu(\beta, h_c^a(\beta) + \Delta) < \mathsf{F}(\beta, h_c^a(\beta) + \Delta)$$

$$\leq \mathsf{F}(0, \Delta) - \frac{\beta^2}{2}(1 - C\beta^2)(\partial_\Delta \mathsf{F}(0, \Delta))^2.$$

The formulation of this theorem has been chosen in order to give a global vision on the *$\mu$ versus $\mathsf{F}$ bounds*, but the novel statement is just the first inequality. The two inequalities in the second line of (2.26) have, in fact, been proven, in [21], Appendix B and [30], Theorem 2.6, respectively.

REMARK 2.5. In the proof of Theorem 2.3, we actually prove the following, more explicit, bound: for every $\varepsilon \in (0, 1)$, there exist positive constants $\beta_0(\varepsilon)$ and $\Delta_0(\varepsilon)$ (we refer to Corollary 5.4 for explicit expressions of these constants) such that

$$(2.27) \qquad \mathsf{F}(\beta, h_c^a(\beta) + \Delta) \geq \mathsf{F}(0, \Delta) - (1 + \varepsilon)\frac{\beta^2}{2}(\partial_\Delta \mathsf{F}(0, \Delta))^2$$

for every $\beta \in [0, \beta_0(\varepsilon)]$ and $\Delta \in [0, \Delta_0(\varepsilon)]$. This estimate must, of course, be matched with the rightmost inequality in (2.26).

We now give two further results that can be considered corollaries of Theorem 2.4 and that show how far the method we are using can be pushed.

The first result gives sharp estimates on the finite-volume free energy at the critical point. First, note that for the annealed system,

$$(2.28) \qquad \frac{1}{N} \log \mathbb{E} Z_{N,\omega}(\beta, h_c^a(\beta)) = \frac{1}{N} \log \mathbf{P}(N \in \tau) \overset{N \to \infty}{\sim} -(1-\alpha) \frac{\log N}{N},$$

where we have used (2.3) in the last step. A heuristic weak-disorder expansion established in [12] suggests that for $\beta$ small, $(1/N)\mathbb{E} \log Z_{N,\omega}(\beta, h_c(\beta))$ has the same behavior for large $N$. Indeed, we can prove the following:

PROPOSITION 2.6.   *If (2.11) holds, then, for $\beta$ sufficiently small,*

$$(2.29) \qquad \left| \frac{1}{N} \mathbb{E} \log Z_{N,\omega}(\beta, h_c(\beta)) - \frac{1}{N} \log \mathbf{P}(N \in \tau) \right| = O\!\left(\frac{1}{N}\right).$$

The second and final result is as follows:

PROPOSITION 2.7.   *Let $n(N) \in \mathbb{N}$ be such that $\lim_{N \to \infty} n(N) = \lim_{N \to \infty}(N - n(N)) = +\infty$. If (2.11) holds, then, for $\beta$ sufficiently small,*

$$(2.30) \qquad \limsup_{h \searrow h_c(\beta)} \lim_{N \to \infty} \frac{\mathbb{E}[\mathbf{P}_{N,\omega}(n(N) \in \tau)^2]}{(\mathbb{E}\mathbf{P}_{N,\omega}(n(N) \in \tau))^2} < \infty.$$

As we are going to explain in a moment, the latter result establishes the *absence of multifractality of the order parameter*. First, the conditions on $n(N)$ simply guarantee that we are looking at a site *in the bulk*, that is, a site whose distance from the boundaries of the system diverges in the thermodynamic limit. To make this clear, think of the case $n(N) = \lfloor N/2 \rfloor$ ($\lfloor \cdot \rfloor$ is the integer part of $\cdot$). Since

$$(2.31) \qquad \lim_{N \to \infty} \mathbb{E}\mathbf{P}_{N,\omega}(\lfloor N/2 \rfloor \in \tau) = \partial_h \mathtt{F}(\beta, h)$$

(see Section 4), it is clear that both numerator and denominator in (2.30) vanish in the limit $N \to \infty$ followed by $h \searrow h_c(\beta)$ and the random variable $\mathbf{P}_{N,\omega}(\lfloor N/2 \rfloor \in \tau)$ tends to zero in probability. This is, however, not enough to guarantee that the ratio in (2.30) remains finite. For instance, in cases where (2.11) does not hold [in particular, when $\alpha = 1/2$ and $L(\cdot)$ is asymptotically constant or $\alpha = 3/4$] recent numerical simulations [26] seem to indicate that the numerator and denominator of (2.30) behave, respectively, like $(h - h_c(\beta))^{y_2}$ and $(h - h_c(\beta))^{y_1}$ with $y_1 > y_2 > 0$. This fact, which is referred to as *multifractality of the order parameter at criticality* [26], would imply, in particular, that the limsup in (2.30) is infinite when disorder is relevant.

2.7. *Interpolation method and the 2-replica homogeneous model.* The basic idea developed in [30] is to introduce the modified free energy

$$(2.32) \quad \mathcal{F}_N(t, \lambda) := \frac{1}{N} \mathbb{E} \log \mathbf{E}^{\otimes 2}[\exp(\mathcal{H}_{N,\omega,t,\lambda,\Delta}(\tau, \tau')); N \in \tau \cap \tau'],$$

where

$$(2.33) \quad \begin{aligned} &\mathcal{H}_{N,\omega,t,\lambda,\Delta}(\tau, \tau') \\ &:= \sum_{n=1}^{N} (\sqrt{t}\beta\omega_n + \Delta - t(\beta^2/2))(\mathbf{1}_{n\in\tau} + \mathbf{1}_{n\in\tau'}) + \lambda\beta^2 \sum_{n=1}^{N} \mathbf{1}_{n\in\tau\cap\tau'}, \end{aligned}$$

$\Delta \geq 0$, $0 \leq t \leq 1$ and $\lambda \geq 0$. It is clear that

$$(2.34) \quad \mathcal{F}_N(0, 0) = 2\mathtt{F}_N(0, \Delta)$$

and

$$(2.35) \quad \mathcal{F}_N(1, 0) = 2\mathtt{F}_N(\beta, h_c^a(\beta) + \Delta).$$

Moreover, via Gaussian integration by parts, it is proven in [30] that

$$(2.36) \quad \frac{\mathrm{d}}{\mathrm{d}t}\mathcal{F}_N(t, \lambda) \leq \frac{\mathrm{d}}{\mathrm{d}\lambda}\mathcal{F}_N(t, \lambda)$$

and

$$(2.37) \quad \begin{aligned} \mathtt{F}_N(\beta, h_c^a(\beta) + \Delta) - \mathtt{F}_N(0, \Delta) &= \frac{1}{2}(\mathcal{F}_N(1, 0) - \mathcal{F}_N(0, 0)) \\ &\geq -\frac{e-1}{2}[\mathcal{F}_N(0, 2) - \mathcal{F}_N(0, 0)]. \end{aligned}$$

Since the second line does not involve the disorder $\omega$, (2.37) is a powerful tool to obtain explicit free energy lower bounds (it is through (2.37) that (2.13) is proved in [30]).

At the heart of the technical arguments in this work is, therefore, a homogeneous 2-replica pinning model which may be of interest in its own right and which we discuss in the remainder of this section. It is indexed by a nonnegative parameter $b$ and by a real coupling parameter $\lambda$, and it is built as follows: recall that we are assuming $\sum_{n\in\mathbb{N}} K(n) = 1$ and set

$$(2.38) \quad K_b(n) := c(b)\exp(-bn)K(n) = c(b)L(n)\frac{\exp(-bn)}{n^{1+\alpha}}$$

with

$$(2.39) \quad c(b)^{-1} := \sum_{n\in\mathbb{N}} K(n)\exp(-bn)$$

so that $c(0) = 1$ and the $K_b(\cdot)$-renewal is persistent. Now, take two independent copies $\tau$ and $\tau'$ of a $K_b(\cdot)$-renewal. As we have already seen,

$$(2.40) \qquad \mathbf{P}_b^{\otimes 2}(n \in \tau \cap \tau') = (\mathbf{P}_b(n \in \tau))^2$$

and through (2.2), the inter-arrival law of $\tau \cap \tau'$ that we denote by $\mathbb{K}_b(\cdot)$ can therefore be computed.

REMARK 2.8. As we have already mentioned, while the renewal function is explicit, $\mathbb{K}_b(\cdot)$ is not, and this leads to substantial difficulties. In particular, while it is not difficult to see that $\limsup_n n^{-1} \log \mathbb{K}_b(n) < 0$ for $b > 0$, identifying the correct rate of convergence or the precise asymptotic behavior of $\mathbb{K}_b(\cdot)$ is harder and results are a priori counterintuitive (see Section 5.1). In reality, when $b > 0$, the sharp decay rate of $\mathbb{K}_b(\cdot)$ is not crucial for us and what turns out to be central (e.g., in the proof of Proposition 5.3) is correctly taking into account the contribution of $\mathbb{K}_b(n)$ for $n = O(1/b)$; one can actually show that in this regime, $\mathbb{K}_b(n)$ significantly differs from what one extrapolates from the tail behavior.

The 2-replica homogeneous model is defined in strict analogy with the 1-replica homogeneous model (i.e., the annealed model) of Section 2.3, with the notational change of $\lambda$ in place of $h$ and the more substantial change of considering $\mathbb{K}_b(\cdot)$-renewals instead of $K(\cdot)$-renewals. We will actually be interested in the free energy of the model, namely, in

$$(2.41) \quad \mathtt{B}(b, \lambda) := \lim_{N \to \infty} \frac{1}{N} \log \mathbf{E}_b^{\otimes 2} \left[ \exp\left( \lambda \sum_{n=1}^{N} \mathbf{1}_{n \in \tau \cap \tau'} \right); N \in \tau \cap \tau' \right].$$

In [30], the term on the right-hand side of (2.41) has been bounded above by using the Hölder inequality, thereby obtaining (2.13). One of the purposes of our work is to sharply evaluate $\mathtt{B}(b, \lambda)$.

Note that, setting $b = 0$ in (2.41), one gets precisely (2.16) and that, thanks to (2.7) and (2.38) (see also Remark 3.1 below),

$$(2.42) \qquad \lim_{N \to \infty} (\mathcal{F}_N(0, \lambda) - \mathcal{F}_N(0, 0)) = \mathtt{B}(\mathtt{F}(0, \Delta), \beta^2 \lambda).$$

The existence of the limit in (2.41) falls in the realm of the theory of the standard homogeneous pinning model outlined above, with the important difference that if $b > 0$, then $\mathbb{K}_b(\cdot)$ decays exponentially and $\mathtt{B}(b, \lambda)$ can then be negative; see Remark 2.8. In fact, once again, the crucial point is to decide whether or not the equation

$$(2.43) \qquad \sum_{n \in \mathbb{N}} \mathbb{K}_b(n) \exp(-\mathtt{B}n) = e^{-\lambda},$$

has a solution. For the sake of conciseness, we will not consider this problem in full generality. The relevant case for us is the one in which (2.11) holds and $\lambda > 0$ [which immediately implies that $\mathtt{B}(b, \lambda) > 0$]. In particular, we will be interested in the asymptotic behavior of $\mathtt{B}(b, \lambda)$ for both $\lambda$ and $b$ small and we will show that in this regime,

$$
(2.44) \qquad \mathtt{B}(b, \lambda) \sim \frac{\lambda}{(\sum_{n \in \mathbb{N}} n K_b(n))^2},
$$

which, incidentally, is just what one would obtain by naively expanding to first order the exponential in (2.41) for $\lambda$ small and then applying the renewal theorem to take the $N \to \infty$ limit (thus performing an a priori unjustified exchange of limits). For a more precise and refined statement of (2.44), see Corollary 5.4 below. The 2-replica model is treated in detail in Section 5.

## 3. Interpolation procedure and proof of the main results.

3.1. *Proof of Theorem 2.4.* From now on, for ease of notation, we set $\delta_n := \mathbf{1}_{n \in \tau}$. As pointed out immediately after the statement, we just need to prove the first inequality in (2.26) and we start from the identity

$$
(3.1) \qquad
\begin{aligned}
&Z_{N, \omega}(\beta, h_c^a(\beta) + \Delta) \\
&\quad = \mathbf{E}(e^{\Delta \sum_{n=1}^{N} \delta_n} \delta_N) \left\langle \exp\left( \sum_{n=1}^{N} (\beta \omega_n - (\beta^2/2)) \delta_n \right) \right\rangle_{N, \Delta},
\end{aligned}
$$

where $\langle \cdot \rangle_{N, \Delta}$ is the measure

$$
(3.2) \qquad \langle f \rangle_{N, \Delta} := \frac{\mathbf{E}(e^{\Delta \sum_{n=1}^{N} \delta_n} f \delta_N)}{\mathbf{E}(e^{\Delta \sum_{n=1}^{N} \delta_n} \delta_N)},
$$

that is, the Boltzmann–Gibbs measure for the homogeneous system with pinning parameter $\Delta > 0$.

REMARK 3.1. It is well known (see, for instance, [18], Section 2.2) that $\langle \cdot \rangle_{N, \Delta}$ can be equivalently described as follows. Let $\tau$, with law $\mathbf{P}_{\mathtt{F}(0, \Delta)}(\cdot)$, be the positive recurrent renewal with inter-arrival law [maintaining consistency with the notation (2.38)]

$$
(3.3) \qquad K_{\mathtt{F}(0, \Delta)}(n) := e^{-n \mathtt{F}(0, \Delta)} K(n) e^{\Delta}.
$$

The fact that $K_{\mathtt{F}(0, \Delta)}(\cdot)$ thus defined is actually normalized to one follows from (2.7). Then,

$$
(3.4) \qquad \langle \cdot \rangle_{N, \Delta} = \mathbf{E}_{\mathtt{F}(0, \Delta)}(\cdot | N \in \tau).
$$

Going back to (3.1), one finds

$$
\begin{aligned}
(3.5) \quad & \mu_N(\beta, h_c^a(\beta) + \Delta) \\
& = \mathsf{F}_N(0, \Delta) - \frac{1}{N} \log \mathbb{E}\left[\left\langle \exp\left(\sum_{n=1}^N (\beta \omega_n - (\beta^2/2))\delta_n\right)\right\rangle_{N,\Delta}^{-1}\right].
\end{aligned}
$$

This leads us to introduce, for $t \in [0, 1]$,

$$
(3.6) \quad \phi_{N,t} := -\frac{1}{N} \log \mathbb{E}\left[\left\langle \exp\left(\sum_{n=1}^N (\beta \sqrt{t} \omega_n - t(\beta^2/2))\delta_n\right)\right\rangle_{N,\Delta}^{-1}\right].
$$

Note that $\phi_{N,0} = 0$, while $\phi_{N,1}$ coincides, by (3.5), with $\mu_N(\beta, h_c^a(\beta) + \Delta) - \mathsf{F}_N(0, \Delta)$. Using integration by parts with respect to the Gaussian density, we obtain

$$
(3.7) \quad \frac{\mathrm{d}}{\mathrm{d}t} \phi_{N,t} = -\frac{\beta^2}{N} \sum_{n=1}^N \frac{\mathbb{E}[(\langle\!\langle \delta_n \rangle\!\rangle_t)^2 / Z_t]}{\mathbb{E}[1/Z_t]},
$$

where we have used the notation

$$
(3.8) \quad \langle\!\langle \cdot \rangle\!\rangle_t := \frac{1}{Z_t} \left\langle \exp\left(\sum_{n=1}^N (\beta \sqrt{t} \omega_n - t(\beta^2/2))\delta_n\right) \cdot \right\rangle_{\Delta,N}
$$

and

$$
(3.9) \quad Z_t := \left\langle \exp\left(\sum_{n=1}^N (\beta \sqrt{t} \omega_n - t(\beta^2/2))\delta_n\right) \right\rangle_{\Delta,N}.
$$

As a convenient notational shortcut, we set $Q_t := Z_t^{-1}/\mathbb{E}[Z_t^{-1}]$. The basic technical estimate is the following:

LEMMA 3.2. *For every choice of $t \in [0, 1]$, $\beta$, $h$ and $N$, we have*

$$
(3.10) \quad \beta^2 \sum_{n=1}^N \mathbb{E}[\langle\!\langle \delta_n \rangle\!\rangle_t^2 Q_t] \le \log\left\langle \exp\left(8\beta^2 \sum_{n=1}^N \delta_n \delta_n'\right)\right\rangle_{N,\Delta}^{\otimes 2}.
$$

By combining the lemma with (3.6) and (3.7), we directly obtain

$$
(3.11) \quad \mu(\beta, h_c^a(\beta) + \Delta) \ge \mathsf{F}(0, \Delta) - \lim_{N \to \infty} \frac{1}{N} \log\left\langle \exp\left(8\beta^2 \sum_{n=1}^N \delta_n \delta_n'\right)\right\rangle_{N,\Delta}^{\otimes 2}.
$$

In view of Remark 3.1 and of the positive recurrence of $\tau$ under $\mathbf{P}_{\mathsf{F}(0,\Delta)}(\cdot)$, the limit in the right-hand side can be written as

$$
\begin{aligned}
(3.12) \quad & \lim_{N \to \infty} \frac{1}{N} \log \mathbf{E}_{\mathsf{F}(0,\Delta)}^{\otimes 2}\left[\exp\left(8\beta^2 \sum_{n=1}^N \mathbf{1}_{\{n \in \tau \cap \tau'\}}\right); N \in \tau \cap \tau'\right] \\
& = \mathsf{B}(\mathsf{F}(0, \Delta), 8\beta^2).
\end{aligned}
$$

From Corollary 5.4 and the observations of Remark 5.2, one concludes that, for $\beta$ and $\Delta$ sufficiently small, the limit (3.12) is smaller than $9\beta^2(\partial_\Delta F(0,\Delta))^2$.

This establishes the first inequality in (2.26) and therefore Theorem 2.4 is proven.

PROOF OF LEMMA 3.2. The proof has much in common with that of Theorem 3.1 in [29], which is attributed by the author of that article to R. Latala. Let us set, for $\lambda \in \mathbb{R}$,

$$(3.13) \qquad \psi_{N,t} := \log \mathbb{E}[Q_t \langle\!\langle \exp(\lambda\beta^2 \underline{\delta} \cdot \underline{\delta}') \rangle\!\rangle^{\otimes 2}],$$

where, of course, $\underline{\delta} \cdot \underline{\delta}' := \sum_{n=1}^N \delta_n \delta'_n$. We will make use of the identity

$$
\begin{aligned}
(3.14) \qquad \frac{\mathrm{d}}{\mathrm{d}t}\psi_{N,t} = {}& (\mathbb{E}[Q_t\langle\!\langle \exp(\lambda\beta^2 \underline{\delta} \cdot \underline{\delta}')\rangle\!\rangle_t^{\otimes 2}])^{-1} \\
&\times \Bigg\{ 6\beta^2 \mathbb{E}[Q_t\langle\!\langle \exp(\lambda\beta^2 \underline{\delta}\cdot\underline{\delta}')\rangle\!\rangle_t^{\otimes 2}\langle\!\langle \underline{\delta}\cdot\underline{\delta}'\rangle\!\rangle_t^{\otimes 2}] \\
&\quad - 6\beta^2 \mathbb{E}\left[ Q_t \sum_{n=1}^N \langle\!\langle \delta_n \exp(\lambda\beta^2 \underline{\delta}\cdot\underline{\delta}')\rangle\!\rangle_t^{\otimes 2}\langle\!\langle \delta_n\rangle\!\rangle_t \right] \\
&\quad + \beta^2 \mathbb{E}[Q_t\langle\!\langle \underline{\delta}\cdot\underline{\delta}' \exp(\lambda\beta^2 \underline{\delta}\cdot\underline{\delta}')\rangle\!\rangle_t^{\otimes 2}] \\
&\quad - \beta^2 \mathbb{E}\left[ Q_t \sum_{n=1}^N (\langle\!\langle \delta_n\rangle\!\rangle_t)^2 \right] \mathbb{E}[Q_t\langle\!\langle \exp(\lambda\beta^2 \underline{\delta}\cdot\underline{\delta}')\rangle\!\rangle_t^{\otimes 2}] \Bigg\}
\end{aligned}
$$

that holds for every $\lambda$ (see below for the steps leading to this identity). Using the fact that $\mathbb{E}(X\exp(aX)) \geq \mathbb{E}(X)\mathbb{E}(\exp(aX))$ for any random variable $X$ [which follows from the monotonicity of $\mathbb{E}(X\exp(aX))/\mathbb{E}(\exp(aX))$ with respect to $a$], one finds that

$$
\begin{aligned}
(3.15) \qquad & \mathbb{E}[Q_t\langle\!\langle \exp(\lambda\beta^2 \underline{\delta}\cdot\underline{\delta}')\rangle\!\rangle_t^{\otimes 2}\langle\!\langle \underline{\delta}\cdot\underline{\delta}'\rangle\!\rangle_t^{\otimes 2}] \\
&\qquad \leq \mathbb{E}[Q_t\langle\!\langle \underline{\delta}\cdot\underline{\delta}' \exp(\lambda\beta^2 \underline{\delta}\cdot\underline{\delta}')\rangle\!\rangle_t^{\otimes 2}].
\end{aligned}
$$

Therefore,

$$
\begin{aligned}
(3.16) \qquad \frac{\mathrm{d}}{\mathrm{d}t}\psi_{N,t} \leq {}& 7\beta^2 (\mathbb{E}[Q_t\langle\!\langle \exp(\lambda\beta^2 \underline{\delta}\cdot\underline{\delta}')\rangle\!\rangle_t^{\otimes 2}])^{-1} \\
&\times \mathbb{E}[Q_t\langle\!\langle \underline{\delta}\cdot\underline{\delta}' \exp(\lambda\beta^2 \underline{\delta}\cdot\underline{\delta}')\rangle\!\rangle_t^{\otimes 2}].
\end{aligned}
$$

We readily see that the previous inequality implies

$$(3.17) \qquad \frac{\mathrm{d}}{\mathrm{d}t} \log \mathbb{E}[Q_t\langle\!\langle \exp(\beta^2(8-7t)\underline{\delta}\cdot\underline{\delta}')\rangle\!\rangle_t^{\otimes 2}] \leq 0$$

so that

$$(3.18) \qquad \log \mathbb{E}[Q_t\langle\!\langle \exp(\beta^2\underline{\delta}\cdot\underline{\delta}')\rangle\!\rangle_t^{\otimes 2}] \leq \log \langle \exp(8\beta^2\underline{\delta}\cdot\underline{\delta}')\rangle_{N,\Delta}^{\otimes 2}$$

for every $t \in [0,1]$. Since Jensen's inequality guarantees that the left-hand side of (3.18) is bounded below by $\beta^2 \mathbb{E}[Q_t \langle\!\langle \underline{\delta} \cdot \underline{\delta}' \rangle\!\rangle_t^{\otimes 2}]$, we have obtained (3.10).

We conclude the proof by giving some details on the computation leading to (3.10). We write

$$\frac{\mathrm{d}}{\mathrm{d}t}\psi_{N,t} = (\mathbb{E}[Q_t \langle\!\langle \exp(\lambda\beta^2 \underline{\delta} \cdot \underline{\delta}') \rangle\!\rangle_t^{\otimes 2}])^{-1} A \tag{3.19}$$

with

$$A := \frac{\mathrm{d}}{\mathrm{d}t} \frac{\mathbb{E}[Z_t^{-3} \langle \exp(\underline{\eta}(t) \cdot (\underline{\delta} + \underline{\delta}') + \lambda\beta^2 \underline{\delta} \cdot \underline{\delta}') \rangle_{N,\Delta}]}{\mathbb{E}[Z_t^{-1}]}, \tag{3.20}$$

where we have use the shorthand notation $\underline{\eta}(t) = (\eta_1(t), \ldots, \eta_N(t))$ with $\eta_n(t) := \beta\sqrt{t}\omega_n - t(\beta^2/2)$. We have, from (3.7),

$$\frac{\mathrm{d}}{\mathrm{d}t}\mathbb{E}[Z_t^{-1}] = \beta^2 \mathbb{E}\left[ Z_t^{-1} \sum_{n=1}^{N} (\langle\!\langle \delta_n \rangle\!\rangle_t)^2 \right]. \tag{3.21}$$

Moreover,

$$\frac{\mathrm{d}}{\mathrm{d}t}\mathbb{E}[Z_t^{-3} \langle \exp(\underline{\eta}(t) \cdot (\underline{\delta} + \underline{\delta}') + \lambda\beta^2 \underline{\delta} \cdot \underline{\delta}') \rangle_{N,\Delta}^{\otimes 2}]$$

$$= -3\mathbb{E}\left[ Z_t^{-4} \langle \exp(\underline{\eta}(t) \cdot (\underline{\delta} + \underline{\delta}') + \lambda\beta^2 \underline{\delta} \cdot \underline{\delta}') \rangle_{N,\Delta}^{\otimes 2} \frac{\mathrm{d}\underline{\eta}(t)}{\mathrm{d}t} \right.$$

$$\left. \cdot \langle \underline{\delta} \exp(\underline{\eta}(t) \cdot \underline{\delta}) \rangle_{N,\Delta} \right] \tag{3.22}$$

$$+ \mathbb{E}\left[ Z_t^{-3} \frac{\mathrm{d}\underline{\eta}(t)}{\mathrm{d}t} \cdot \langle (\underline{\delta} + \underline{\delta}') \exp(\underline{\eta}(t) \cdot (\underline{\delta} + \underline{\delta}') + \lambda\beta^2 \underline{\delta} \cdot \underline{\delta}') \rangle_{N,\Delta}^{\otimes 2} \right]$$

and, by (Gaussian) integration by parts, we see that the right-hand side is equal to

$$6\beta^2 \mathbb{E}\left[ Z_t^{-1} \langle\!\langle \exp(\lambda\beta^2 \underline{\delta} \cdot \underline{\delta}') \rangle\!\rangle_t^{\otimes 2} \sum_{n=1}^{N} (\langle\!\langle \delta_n \rangle\!\rangle_t)^2 \right]$$

$$- 6\beta^2 \mathbb{E}\left[ Z_t^{-1} \sum_{n=1}^{N} \langle\!\langle \delta_n \exp(\lambda\beta^2 \underline{\delta} \cdot \underline{\delta}') \rangle\!\rangle_t^{\otimes 2} \langle\!\langle \delta_n \rangle\!\rangle_t \right] \tag{3.23}$$

$$+ \beta^2 \mathbb{E}[Z_t^{-1} \langle\!\langle \underline{\delta} \cdot \underline{\delta}' \exp(\lambda\beta^2 \underline{\delta} \cdot \underline{\delta}') \rangle\!\rangle_t^{\otimes 2}].$$

Inserting the last expressions into the definition of $A$ [see (3.20)] and, in turn, into (3.19) one obtains (3.14).   $\square$

3.2. *Proof of Theorem 2.3.* By the last inequality in (2.26), we already know that, for $\Delta$ and $\beta$ small enough,

$$(3.24) \qquad \frac{\mathtt{F}(0,\Delta) - \mathtt{F}(\beta, h_c^a(\beta) + \Delta)}{(\beta^2/2)(\partial_\Delta \mathtt{F}(0,\Delta))^2} - 1 \geq -\beta^2 C.$$

To obtain the complementary bound, we recall the bound (2.37) proven in [30], Section 3. In such a bound, the multiplicative factor $e - 1$ appears because we have chosen to let $\lambda$ run from 1 to 2 [and this yields the constant 2 appearing in the term $\mathcal{F}_N(0,2)$]. Letting instead $\lambda$ run up to $M + 1$ ($M$ is going to be chosen to be large) in that proof allows us to replace, in the $N \to \infty$ limit, the last term in the right-hand side of (2.37) by

$$(3.25) \qquad -\frac{(e^{1/M} - 1)}{2} \lim_{N \to \infty} [\mathcal{F}_N(0, 1 + M) - \mathcal{F}_N(0,0)]$$

for every $M > 0$ so that [see (2.42)]

$$(3.26) \qquad \frac{\mathtt{F}(0,\Delta) - \mathtt{F}(\beta, h_c^a(\beta) + \Delta)}{\mathtt{B}(\mathtt{F}(0,\Delta), \beta^2(1 + M))(\exp(1/M) - 1)/2} - 1 \leq 0.$$

By applying Corollary 5.4 and recalling Remark 5.2, we see that for every $\beta > 0$ and every $M$ for $\lambda$ and $\Delta$ sufficiently small, the denominator in the left-hand side of (3.26) is bounded above by

$$(3.27) \qquad \tfrac{1}{2}\beta^2(\partial_\Delta \mathtt{F}(0,\Delta))^2 (\exp(1/M) - 1)(1 + M)(1 + \varepsilon).$$

The proof of (2.15) is thus complete because we can choose $\varepsilon > 0$ and $1/M$ arbitrarily small.

**4. Proof of the corollaries to the main results.** In this section, we are going to prove Propositions 2.6 and 2.7.

PROOF OF PROPOSITION 2.6. We only need to prove that

$$(4.1) \qquad \mathtt{F}_N(\beta, h_c^a(\beta)) \geq \frac{1}{N} \log \mathbf{P}(N \in \tau) + O\!\left(\frac{1}{N}\right)$$

since the complementary bound (without any error term) is simply Jensen's inequality. On the other hand, (2.37) for $\Delta = 0$ gives

$$(4.2) \qquad \mathtt{F}_N(\beta, h_c^a(\beta)) \geq \frac{1}{N} \log \mathbf{P}(N \in \tau) - \frac{e-1}{2N} \log \mathbf{E}^{\otimes 2}(e^{2\beta^2 \underline{\hat{\Sigma}} \cdot \underline{\hat{\Sigma}}'} | N \in \tau^1 \cap \tau^2).$$

From (2.3) and Proposition 5.7, we know that $\tau^1 \cap \tau^2$ is a terminating renewal whose inter-arrival law $\mathbb{K}_0(\cdot)$ satisfies

$$(4.3) \qquad \mathbb{K}_0(n) \overset{n \to \infty}{\sim} \frac{C}{L(n)^2 n^{2-2\alpha}}$$

for some positive constant $C$ which depends on $\alpha$ and $L(\cdot)$. Therefore, for $\lambda$ smaller than $-\log(\sum_{n\in\mathbb{N}}\mathbb{K}_0(n))$ (and, in particular, for $\lambda = 0$), we have

$$(4.4) \quad \frac{1}{2N}\log \mathbf{E}^{\otimes 2}(e^{\lambda\underline{\delta}\cdot\underline{\delta}'}; N \in \tau^1 \cap \tau^2) = -(1-\alpha)\frac{\log N}{N} + O\left(\frac{1}{N}\right).$$

This result follows by observing that we are pinning a terminating renewal with regularly varying inter-arrival distribution: the condition on $\lambda$ is precisely given to ensure that such a pinning is not sufficient to localize the model. As shown in, for example, [18], Theorem 2.2, in such a regime, the partition function behaves, to leading order, like the inter-arrival law $\mathbb{K}_0(N)$, up to an explicit (in this case $\lambda$-dependent) multiplicative constant. Equation (2.29) therefore follows, provided that $\beta$ is small enough.  $\square$

PROOF OF PROPOSITION 2.7.   As usual, we set $\Delta = h - h_c(\beta) = h - h_c^a(\beta)$ and we begin by observing that

$$(4.5) \quad \lim_{N\to\infty}\mathbb{E}\mathbf{P}_{N,\omega}(n(N) \in \tau) = \partial_h\mathtt{F}(\beta, h) \overset{\Delta\to 0}{\sim} \partial_\Delta\mathtt{F}(0, \Delta).$$

The first equality and the existence of the limit follows from the exponential decay of correlations in the localized phase [21], while the asymptotic equality for $\Delta \to 0$ follows from Theorem 2.4 together with the convexity of $\mathtt{F}(0, \cdot)$.

The claim of the theorem follows once we show that for every $\beta$ (sufficiently small), there exists a constant $C > 0$ such that, say for $0 < \Delta < 1$,

$$(4.6) \quad \lim_{N\to\infty}\mathbb{E}[(\mathbf{P}_{N,\omega}(n(N) \in \tau))^2] \leq C(\partial_\Delta\mathtt{F}(0, \Delta))^2.$$

To prove this, we first observe that, again thanks to the exponential decay of correlations,

$$
\begin{aligned}
(4.7) \quad & \lim_{N\to\infty}\mathbb{E}[(\mathbf{P}_{N,\omega}(n(N) \in \tau))^2] \\
& = \lim_{N\to\infty}\frac{1}{N}\sum_{n=1}^{N}\mathbb{E}[(\mathbf{P}_{N,\omega}(n \in \tau))^2] \\
& = \partial_\lambda\frac{1}{N}\mathbb{E}\log\mathbf{E}^{\otimes 2}(e^{\sum_n(\beta\omega_n+h)(\delta_n+\delta'_n)+\lambda\underline{\delta}\cdot\underline{\delta}'})|_{\lambda=0}.
\end{aligned}
$$

Using Jensen's inequality, the right-hand side is bounded above, for every $\lambda > 0$, by

$$
\begin{aligned}
(4.8) \quad & \frac{1}{\lambda}\left[\frac{1}{N}\mathbb{E}\log\mathbf{E}^{\otimes 2}(e^{\sum_n(\beta\omega_n+h)(\delta_n+\delta'_n)+\lambda\underline{\delta}\cdot\underline{\delta}'}) - 2\mathtt{F}_N(\beta, h)\right] \\
& = \frac{1}{\lambda}\left[\mathcal{F}_N\left(1, \frac{\lambda}{\beta^2}\right) - 2\mathtt{F}_N(\beta, h)\right].
\end{aligned}
$$

Using (2.42) and (2.36), one therefore finds that

$$
\begin{aligned}
(4.9) \quad & \lim_{N \to \infty} \mathbb{E}[(\mathbf{P}_{N,\omega}(n(N) \in \tau))^2] \\
& \leq \frac{1}{\lambda}[2(\mathbf{F}(0,\Delta) - \mathbf{F}(\beta,h)) + \mathbf{B}(\mathbf{F}(0,\Delta),\lambda + \beta^2)].
\end{aligned}
$$

Now choose, for example, $\lambda = \beta^2$ and apply Theorem 2.4, Corollary 5.4 and (5.13) to obtain (4.6). This concludes the proof of Proposition 2.7. $\quad\square$

## 5. On intersection two independent renewals.

In this section, we study some properties of renewals that are obtained as intersections of two independent copies of a given $K_b(\cdot)$-renewal. The results we obtain may be of independent interest and, for this reason, this section is somewhat independent of the rest of the work.

We are going to consider renewal processes $\tau = \{\tau_j\}_{j=0,1,\ldots}$ with inter-arrival law supported by $\mathbb{N} = \{1,2,\ldots\}$. The renewal function of $\tau$ computed in $n \in \mathbb{N} \cup \{0\}$ is, by definition, $\mathbf{P}(n \in \tau)$ and the renewal function is related to the inter-arrival distribution [i.e., to the function $n \mapsto \mathbf{P}(\tau_1 = n)$] by the recurrence scheme (2.2) for $n = 0,1,\ldots$. In what follows, a renewal process such that $\mathbf{P}(\tau_1 = n) = F(n)$ for every $n$ is called $F(\cdot)$-*renewal* and $n \mapsto \mathbf{P}(n \in \tau)$ is the corresponding $F(\cdot)$-*renewal function*.

For $b \geq 0$ and $n \in \mathbb{N}$, let $K_b(\cdot)$ be defined as in (2.38), where $K(\cdot) = K_0(\cdot)$ satisfies (2.1) and we assume (2.11). We define

$$
(5.1) \qquad \mathbb{K}_b(n) := \mathbf{P}_b^{\otimes 2}((\tau \cap \tau')_1 = n),
$$

the return distribution of the intersection of two independent $K_b(\cdot)$-renewals, and

$$
(5.2) \qquad u_b(n) := \mathbf{P}_b(n \in \tau),
$$

the renewal function of a single $K_b(\cdot)$-renewal. If $b > 0$, since the two renewals are positive recurrent, the intersection is positive recurrent too. If $b = 0$, the two renewals are null-recurrent, but the intersection is terminating, as was discussed just before (2.16).

REMARK 5.1. Note that

$$
\begin{aligned}
(5.3) \quad u_b(\infty) &:= \lim_{n \to \infty} u_b(n) = \frac{1}{\sum_n c(b)L(n)n^{-\alpha}\exp(-bn)} \\
& \stackrel{b \searrow 0}{\sim} \frac{b^{1-\alpha}}{L(1/b)\Gamma(1-\alpha)},
\end{aligned}
$$

where the equality is a consequence of the renewal theorem. So, in particular,

$$
(5.4) \qquad u_b^2(\infty) = o(b).
$$

The latter statement is also true if $\alpha = 1/2$ and $L(\cdot)$ diverges at infinity and we remark that (2.11) does indeed imply that $\lim_n L(n) = \infty$. Let us also observe that for the normalization constant $c(b)$ in (2.38), we have $c(b) - 1 \sim b^\alpha L(1/b)\Gamma(1-\alpha)/\alpha$ so that $1/(c(b)-1) \sim \alpha u_b(\infty)/b$.

In most situations, one has a good grip on the renewal process if the inter-arrival law is known in detail. In our case, $\mathbb{K}_b(\cdot)$ has been introduced only indirectly and it is natural to try to characterize it as precisely as possible. What is instead characterized in a straightforward way is the $\mathbb{K}_b(\cdot)$-renewal function, which we denote by $\mathbb{U}_b(\cdot)$ since, as we have already mentioned, $\mathbb{U}_b(n) = u_b^2(n)$.

Let us take this opportunity to point out the following identity, which is a direct consequence of (2.2):

$$(5.5) \qquad \widehat{u}_b(z) = \frac{1}{1 - \widehat{K}_b(z)},$$

where we have used the notation $\widehat{f}(z) := \sum_{n=0}^\infty f(n)z^n$, with $z$ in a centered ball of $\mathbb{C}$. Equation (5.5) holds for $|z|$ within the radius of convergence of the two series appearing in the expression (and in any case for $|z| < 1$).

Of course, in (5.5), we can replace $u_b$ with $\mathbb{U}_b$ and $K_b$ with $\mathbb{K}_b$. Looking at the problem this way, retrieving $\mathbb{K}_b(\cdot)$ from $\mathbb{U}_b(\cdot)$ is an inverse $z$-transform problem. We will attack this question in some detail in Section 5.1; the main problem we want to tackle here is computing, for both $\lambda$ and $b$ positive, the limit as $N$ tends to infinity of

$$(5.6) \qquad \frac{1}{N}\log \mathbf{E}_b^{\otimes 2}\left[\exp\left(\lambda \sum_{n=1}^N \mathbf{1}_{n \in \tau \cap \tau'}\right); N \in \tau \cap \tau'\right],$$

where $\tau$ and $\tau'$ are independent copies of a $K_b(\cdot)$-renewal. One can show the existence of this limit and give an expression for its value by applying the procedure detailed in, for example, [18], Chapter 1: note that we are, in fact, just computing the free energy of the homogeneous pinning model, based on the $\mathbb{K}_b(\cdot)$-renewal, with pinning interaction $\lambda$. What we obtain by applying such a procedure is that the limit as $N \to \infty$ of the expression in (5.6) is the unique solution $\mathtt{B} := \mathtt{B}(b, \lambda) > 0$ of

$$(5.7) \qquad \left(\sum_{n=1}^\infty \mathbb{K}_b(n)\exp(-\mathtt{B}n) = \right)\widehat{\mathbb{K}}_b(\exp(-\mathtt{B})) = e^{-\lambda}.$$

Note that existence of the solution is an immediate consequence of the recurrent character of the $\mathbb{K}_b(\cdot)$-renewal $(b > 0)$ and monotonicity. Our aim is to find sharp estimates on $\mathtt{B}(b, \lambda)$ as $b \searrow 0$. Note, also, that if $b = 0$, one can solve the problem in (5.7) only if $\lambda \geq -\log(1 - \widehat{\mathbb{K}}_0(1))$ $[\widehat{\mathbb{K}}_0(1) < 1$ since the $\mathbb{K}_0(\cdot)$-renewal is transient].

In order to achieve our goals, we make some preliminary observations. Note that if we set $\mathbb{D}_b(\cdot) := \mathbb{U}_b(\cdot) - \mathbb{U}_b(\infty)$, by exploiting the basic renewal equation (5.5) applied to the $\mathbb{K}_b(\cdot)$-renewal, we find that

$$(5.8) \qquad \widehat{\mathbb{D}}_b(z) = \frac{1}{1 - \widehat{\mathbb{K}}_b(z)} - \frac{u_b^2(\infty)}{1 - z}.$$

By (5.8), (5.7) is equivalent to

$$(5.9) \qquad \frac{1 - \exp(-\mathtt{B})}{(1 - \exp(-\mathtt{B}))\widehat{\mathbb{D}}_b(\exp(-\mathtt{B})) + u_b^2(\infty)} = 1 - e^{-\lambda}.$$

In view of the asymptotic limit $b \searrow 0$, we make the change of variable

$$(5.10) \qquad x := \frac{1 - \exp(-\mathtt{B})}{u_b^2(\infty)}$$

so that, from (5.9), we obtain

$$(5.11) \qquad x = \frac{1 - e^{-\lambda}}{1 - (1 - e^{-\lambda})\widehat{\mathbb{D}}_b(\exp(-\mathtt{B}))}.$$

REMARK 5.2. An observation that is very relevant to our applications, but not to this section in itself, is that since $\mathtt{F}(h) := \mathtt{F}(0, h)$ is a solution to (2.7) for $h > 0$, it follows that

$$(5.12) \qquad \partial_h \mathtt{F}(0, h) \sum_n \frac{L(n)}{n^\alpha} \exp(-\mathtt{F}(0, h)n) = \exp(-h),$$

which, recalling (5.3) and the fact that $c(\mathtt{F}(0, h)) = e^h$ [see (2.39)], gives, for every $h > 0$,

$$(5.13) \qquad \partial_h \mathtt{F}(0, h) = u_{\mathtt{F}(0, h)}(\infty).$$

For our purposes, the main result of this section is the following:

PROPOSITION 5.3. *For every $c > 0$,*

$$(5.14) \qquad D(c) := \sup_{\mathtt{B} \geq 0, b \in (0, c)} \widehat{\mathbb{D}}_b(\exp(-\mathtt{B})) < \infty.$$

The key to the proof of Proposition 5.3 is controlling $\widehat{\mathbb{D}}_b(\exp(-\mathtt{B}))$ when both $b$ and $\mathtt{B}/b$ are small. Let us also point out that $D(c) > 1$ [see (5.21) and (5.22) below]. But let us first look at an important consequence of Proposition 5.3. We need the auxiliary quantities

$$(5.15) \qquad \lambda_0 := -\log(1 - (2D(c))^{-1})$$

and

$$(5.16) \qquad c_1 := \sup_{\lambda \in (0, \lambda_0]} \frac{1}{2} \left| \frac{\mathrm{d}^2}{\mathrm{d}\lambda^2} \left( \frac{1 - \exp(-\lambda)}{1 - (1 - \exp(-\lambda))D(c)} \right) \right|.$$

Moreover, given $\varepsilon \in (0, 1)$, we set

$$(5.17) \quad b_0(\varepsilon) := \min(c, \inf\{b > 0 : \exp(2(\lambda_0 + c_1\lambda_0^2)u_b^2(\infty)) - 1 \ge \varepsilon\}).$$

Note that $b_0(\varepsilon) > 0$ since $u_b(\infty)$ vanishes as $b \to 0$.

COROLLARY 5.4. *Let us fix $c > 0$. For every $\varepsilon \in (0, 1)$, there exists $b_0(\varepsilon) \in (0, c]$ such that*

$$(5.18) \qquad \mathtt{B}(b, \lambda) \le (1 + \varepsilon)u_b^2(\infty)(\lambda + c_1\lambda^2)$$

*for every $\lambda \le \lambda_0$ and every $b \le b_0(\varepsilon)$.*

PROOF. Going back to (5.11), we see that, by the choice of $\lambda_0$ and $c_1$, we have, for $\lambda \le \lambda_0$,

$$(5.19) \qquad x \le \lambda + c_1\lambda^2$$

and the statement follows from the definition of $x$ [see (5.10)] since, at this point, it is clear that for $B$ sufficiently small, $1 - \exp(-B) \ge B/(1 + \varepsilon)$. The choice we have made for $b_0(\varepsilon)$ guarantees this for $B = \mathtt{B}(b, \lambda)$ since if $1 - \exp(-B) \le \delta$, with $\delta \ge 0$ such that $\exp(2\delta) - 1 \le \varepsilon \in (0, 1)$, then $1 - \exp(-B) \ge B/(1 + \varepsilon)$.

For completeness. we point out that it is immediate to obtain the lower bound

$$(5.20) \qquad \mathtt{B}(b, \lambda) \ge u_b^2(\infty)\lambda$$

that holds in full generality. This is easily obtained by applying Jensen's inequality to (5.6). □

We now turn to Proposition 5.3, but. before starting the proof, we point out that a byproduct of the proof is, in particular, the sharper estimate

$$(5.21) \qquad \lim_{\substack{b \searrow 0 \\ B = o(b)}} \widehat{\mathbb{D}}_b(\exp(-B)) = \widehat{\mathbb{D}}_0(1),$$

where $\widehat{\mathbb{D}}_0(1) := \lim_{B \searrow 0} \widehat{\mathbb{D}}_0(\exp(-B))$, that is,

$$(5.22) \qquad \widehat{\mathbb{D}}_0(1) = \frac{1}{1 - \widehat{\mathbb{K}}_0(1)} = \sum_{n=0}^{\infty} u_0^2(n)$$

and, recalling (2.3) and $u_0(0) = 1$, $\widehat{\mathbb{D}}_0(1) \in (1, \infty)$. Also, note that $\widehat{\mathbb{D}}_0(1)$ can be simply interpreted as the expected size of the $\mathbb{K}_0(\cdot)$-renewal set.

Let us take this opportunity to recall a result that follows directly from the main result (in the discrete setting) of [27]: for every $b > 0$, there exist two positive constants $C_1$ and $C_2$ such that

$$(5.23) \qquad |u_b(n) - u_b(\infty)| \leq C_1 \exp(-C_2 n).$$

The dependence of $C_1$ and $C_2$ on $b$ is rather explicit and one readily sees that one can choose these constants to depend continuously on $b$ so that (5.23) holds also uniformly in $b \in [b_0, b_1]$, for $0 < b_0 < b_1 < \infty$, with $C_1$ and $C_2$ replaced, respectively, by the maximum and the minimum of the same quantities for $b$ ranging in the allowed interval. On the other hand, by using, for instance, the examples in [19], Section 4, one directly sees that there is no choice of $C_1$ and $C_2$ (two constants, not depending on $b$) such that (5.23) holds uniformly in $b \in [b_0, \infty)$. While the estimate (5.23) is relevant for our proof, neither the coupling techniques in [27] nor the precise tail estimates that one obtains by, for example, exploiting [11] (see [19]) are sufficient to control $\widehat{\mathbb{D}}_b(\exp(-B))$ for $b$ small and $B = o(b)$. And the latter regime is the core of the proof, both because it is there that the difficulty lies and because a posteriori in the main application, Corollary 5.4, the free energy $\mathtt{B}(b, \lambda)$ turns out to be $O(u_b^2(\infty))$, whicih is much smaller than $b$.

PROOF OF PROPOSITION 5.3. Let us start by writing the equality

$$\widehat{\mathbb{D}}_b(\exp(-B)) = \sum_{n=0}^{\infty} (u_b(n) - u_b(\infty))^2 \exp(-Bn)$$

$$+ 2 u_b(\infty) \sum_{n=0}^{\infty} (u_b(n) - u_b(\infty)) \exp(-Bn)$$

$$(5.24) \qquad \leq \sum_{n=0}^{\infty} (u_b(n) - u_b(\infty))^2$$

$$+ 2 u_b(\infty) \sum_{n=0}^{\infty} (u_b(n) - u_b(\infty)) \exp(-Bn)$$

$$=: T_1(b) + T_2(b, B).$$

We now claim that it is sufficient to show that

$$(5.25) \qquad \limsup_{b \searrow 0} T_1(b) < \infty$$

and that, given any positive function $B_0(\cdot)$ such that $\lim_{b \searrow 0} B_0(b)/b = 0$, we have

$$(5.26) \qquad \limsup_{b \searrow 0} \sup_{B \in [0, B_0(b)]} |T_2(b, B)| < \infty,$$

to complete the proof of the proposition [we will actually directly identify the limit in (5.25) and show that the limit in (5.26) is zero]. This is because, from (5.24), we obtain that, for $B > 0$,

$$
\begin{aligned}
(5.27) \quad \widehat{\mathbb{D}}_b(\exp(-B)) &\le \sum_{n=0}^{\infty} (u_b(n) - u_b(\infty))^2 \\
&+ 2u_b(\infty) \sqrt{\frac{1}{\exp(2B) - 1} \sum_{n=0}^{\infty} (u_b(n) - u_b(\infty))^2}.
\end{aligned}
$$

Therefore, by using (5.25), we see that $\widehat{\mathbb{D}}_b(\exp(-B))$ is bounded if $u_b(\infty)/\sqrt{B}$ is bounded (smaller than 1, say), at least for $b$ small. But, this means that $\widehat{\mathbb{D}}_b(\exp(-B))$ remains bounded for $b$ small as long as $B > u_b^2(\infty) = o(b)$. And the case $B \le u_b^2(\infty) =: B_0(b)$ is covered precisely by (5.25) and (5.26). This takes care of (5.14) for $b$ sufficiently small, say $b < b_0$.

If, instead, $b \in [b_0, c)$, the bound can be obtained in a rougher way, namely, by observing that $\widehat{\mathbb{D}}_b(\exp(-B)) \le \sum_n |u_b^2(n) - u_b^2(\infty)|$ and using the (uniform) exponential decay of $u_b(n) - u_b(\infty)$ for $b$ in any compact subset of $(0, \infty)$; see (5.23).

Let us therefore go back to (5.25) and (5.26), and let us set $\Delta_b(n) := u_b(n) - u_b(\infty)$.

For $T_2$, we first observe that

$$
(5.28) \quad T_2 = 2u_b(\infty)\widehat{\Delta}_b(\exp(-B))
$$

and that, by (5.3), $T_2$ is bounded if $\widehat{\Delta}_b(\exp(-B)) = O(L(1/b)b^{\alpha-1})$. We will actually show that, uniformly in $B \le B_0(b)$, $\widehat{\Delta}_b(\exp(-B)) = O(b^{-\alpha}/L(1/b))$.

We use the expression

$$
(5.29) \quad \widehat{\Delta}_b(z) = \frac{\widehat{K}_b''(1)}{2(\widehat{K}_b'(1))^2}\left(\frac{F^{(2)}(z)}{F^{(1)}(z)}\right),
$$

where $F^{(1)}(z)$ is the $z$-transform of the proper (that is, not taking the value $\infty$) random variable $\tau^{(1)}$ with law given by $\mathbf{P}(\tau^{(1)} = m) \propto \sum_{n>m} K_b(n)$, $m = 0, 1, \ldots$, and $F^{(2)}(z)$ is the $z$-transform of the proper random variable $\tau^{(2)}$ built by integrating the tail of $\tau^{(1)}$.

A (relatively) straightforward estimate shows that

$$
(5.30) \quad \frac{\widehat{K}_b''(1)}{(\widehat{K}_b'(1))^2} \overset{b \searrow 0}{\sim} \frac{b^{-\alpha}}{L(1/b)}\frac{(1-\alpha)}{\Gamma(1-\alpha)}
$$

and the rest of the argument which will bound $T_2$ is devoted to showing that, in the same limit and uniformly in $B \le B_0(b)$, $F^{(2)}(\exp(-B))/F^{(1)}(\exp(-B)) = O(1)$. Notice that both $F^{(1)}(\cdot)$ and $F^{(2)}(\cdot)$ have a (nonexplicit) dependence on $b$: were they independent of $b$, the result would be immediate (some of the formulae here are given for later use).

Let us start by writing, for $i = 1, 2$,

$$F^{(i)}(z) = \frac{Q_i(z)}{Q_i(1)} \tag{5.31}$$

with

$$Q_1(z) = \sum_{n=0}^{\infty} z^n \sum_{j > n} \frac{L(j)}{j^{1+\alpha}} \exp(-bj) = \sum_{j=1}^{\infty} \frac{L(j)}{j^{1+\alpha}} \exp(-bj) \left( \frac{1 - z^j}{1 - z} \right) \tag{5.32}$$

and

$$\begin{aligned}
Q_2(z) &= \sum_{n=0}^{\infty} z^n \sum_{m > n} \sum_{j > m} \frac{L(j)}{j^{1+\alpha}} \exp(-bj) \\
&= \sum_{j=2}^{\infty} \frac{L(j)}{j^{1+\alpha}} \exp(-bj) \left( \frac{j(1-z) - 1 + z^j}{(1-z)^2} \right).
\end{aligned} \tag{5.33}$$

Moreover,

$$\begin{aligned}
Q_1(1) &= \sum_{j=1}^{\infty} \frac{L(j)}{j^{\alpha}} \exp(-bj) \quad \text{and} \\
Q_2(1) &= \sum_{j=2}^{\infty} \frac{(j-1)L(j)}{2j^{\alpha}} \exp(-bj).
\end{aligned} \tag{5.34}$$

We first note that $|F^{(2)}(z)| \leq 1$ for $|z| \leq 1$. As for $F^{(1)}(z)$, a Taylor expansion of $1 - \exp(-Bj)$ with the bound

$$q := \sup_{x > 0} \left| \frac{1 - e^{-x} - x}{x^2} \right| < \infty \tag{5.35}$$

implies that

$$|Q_1(\exp(-B)) - Q_1(1)| \leq \frac{qB^2}{1 - \exp(-B)} \sum_{j=1}^{\infty} L(j) j^{1-\alpha} \exp(-bj). \tag{5.36}$$

But, $\sum_{j \geq 1} L(j) j^{1-\alpha} \exp(-bj)$ is asymptotically equivalent as $b \searrow 0$ to $b^{\alpha-2} \times \Gamma(2 - \alpha)$, while

$$Q_1(1) \sim L(1/b) b^{\alpha-1} \Gamma(1 - \alpha), \tag{5.37}$$

from which we readily see that there exists some $q_1 > 0$ such that, for $b$ sufficiently small,

$$\left| \frac{Q_1(\exp(-B))}{Q_1(1)} - 1 \right| \leq q_1 \frac{B}{b}. \tag{5.38}$$

This concludes the proof that $\sup_{B \in [0, B_0(b)]} T_2(b, B)$ vanishes as $b \searrow 0$ and, therefore, (5.26) holds.

We are now going to bound $T_1(b)$ by giving a slightly more general argument that, with very little extra effort, is also going to yield (5.21). In what follows, the bound on $T_1(b)$ follows by simply setting $B = 0$.

We start by observing that

$$(5.39) \quad \sum_{n=0}^{\infty} (\Delta_b(n) \exp(-Bn/2))^2 = \int_0^1 |\widehat{\Delta}_b(\exp(-(B/2) + 2\pi i\theta))|^2 \, \mathrm{d}\theta,$$

where we have used Plancherel's formula; for $B = 0$, this expression is, of course, just $T_1(b)$. We now go back to formulae (5.29)–(5.34) and compute. For any fixed $\theta \in (0, 1)$, for $b \searrow 0$ and $B = o(b)$, we have

$$(5.40) \quad Q_1(\exp(-(B/2) + 2\pi i\theta)) \sim \frac{G(\theta)}{1 - \exp(2\pi i\theta)}$$

with

$$(5.41) \quad G(\theta) := \sum_{j=1}^{\infty} \frac{L(j)}{j^{1+\alpha}} (1 - \exp(2\pi i\theta j)).$$

Properties of $G(\cdot)$ are given in Lemma 5.5. Notice, moreover, that $G(t) = \overline{G(1-t)}$ so that the singular behavior at 0 is analogous to the singular behavior at 1. Recalling (5.37), we have

$$(5.42) \quad F^{(1)}(\exp(-(B/2) + 2\pi i\theta)) \sim \frac{G(\theta)}{1 - \exp(2\pi i\theta)} \left( \frac{b^{1-\alpha}}{L(1/b)\Gamma(1-\alpha)} \right).$$

On the other hand, again for every $\theta \in (0, 1)$, we have

$$(5.43) \quad Q_2(\exp(-(B/2) + 2\pi i\theta)) \sim \frac{1}{1 - \exp(2\pi i\theta)} L(1/b) b^{\alpha-1} \Gamma(1-\alpha)$$

and

$$(5.44) \quad Q_2(1) \sim \tfrac{1}{2} b^{\alpha-2} L(1/b)(1-\alpha)\Gamma(1-\alpha).$$

The last two estimates yield

$$(5.45) \quad F^{(2)}(\exp(-(B/2) + 2\pi i\theta)) \sim \frac{2}{1 - \exp(2\pi i\theta)} \frac{b}{(1-\alpha)}.$$

By inserting (5.30), (5.42) and (5.45) into (5.29), we obtain that, for every $\theta \in (0, 1)$,

$$(5.46) \quad \widehat{\Delta}_b(\exp(-(B/2) + 2\pi i\theta)) \sim \frac{1}{G(\theta)}$$

as $b \searrow 0$ with $B = o(b)$. In order to conclude that

$$(5.47) \quad \lim_{\substack{b \to 0 \\ B = o(b)}} \sum_{n=0}^{\infty} (\Delta_b(n) \exp(-Bn/2))^2 = \lim_{b \searrow 0} T_1(b) = \int_0^1 \frac{1}{|G(\theta)|^2} \, \mathrm{d}\theta,$$

we need a domination argument: the bound is detailed in Lemma 5.6. Note that the expression in (5.47) is $\widehat{\mathbb{D}}_0(1)$, which is defined immediately before Proposition 5.3. This is a simple consequence of Plancherel's formula

$$(5.48) \qquad \widehat{\mathbb{D}}_0(1) = \sum_{n \geq 0} u_0(n)^2 = \int_0^1 |\widehat{u}_0(e^{2i\pi\theta})|^2 \, \mathrm{d}\theta$$

and of (5.5).  □

LEMMA 5.5. *We consider the function* $(0,1) \ni \theta \mapsto G(\theta) \in \mathbb{C}$ *defined in (5.41).* $G(\cdot)$ *is continuous and* $|G(\cdot)| > 0$. *Moreover, for every* $\alpha \in (0,1)$, *there exists some* $c_\alpha > 0$ *such that*

$$(5.49) \qquad |G(\theta)| \overset{\theta \searrow 0}{\sim} c_\alpha \theta^\alpha L(1/\theta)$$

*and, therefore,* $1/G(\cdot) \in L^2$ *if* $\alpha \in (0, 1/2)$, *or if* $\alpha = 1/2$ *and* $\int_0^1 (\sqrt{t} L(1/t))^{-2} \, \mathrm{d}t < \infty$.

PROOF. The series defining $G(\cdot)$ is absolutely convergent, so continuity follows by dominated convergence. The positivity of the absolute value follows from the positivity of the real part. The proof of (5.49) follows from a Riemann sum approximation: the asymptotic behavior of the $\Im(G(\cdot))$ is given in [6], 4.3.1a and the real part is treated similarly. In detail,

$$(5.50) \quad G(\theta) \overset{\theta \searrow 0}{\sim} \theta^\alpha L(1/\theta) \left( \int_0^\infty \frac{1 - \cos(2\pi t)}{t^{1+\alpha}} \, \mathrm{d}t + i \int_0^\infty \frac{\sin(2\pi t)}{t^{1+\alpha}} \, \mathrm{d}t \right).$$

More explicitly,

$$(5.51) \qquad \begin{aligned} G(\theta) \overset{\theta \searrow 0}{\sim} {}& \theta^\alpha L(1/\theta)(2\pi)^\alpha \frac{\Gamma(1-\alpha)}{\alpha} \\ & \times \left( \left( \frac{1}{2} \cos(\alpha\pi/2) \right) + i \sin(\alpha\pi/2) \right). \end{aligned}$$  □

LEMMA 5.6. *Assume that (2.11) holds and that* $B = o(b)$ *[as in (5.26)].* *There exist positive constants* $C_1, C_2$ *and* $b_0$ *such that* $C_1 b_0 < 1/2$ *and, for* $0 < b < b_0$,

$$(5.52) \qquad \begin{aligned} &|\widehat{\Delta}_b(\exp(-(B/2) + 2\pi i\theta))| \\ &\leq C_2 \begin{cases} |G(\theta)|^{-1}, & \text{if } \theta \in (C_1 b, 1/2], \\ b^{-\alpha}/L(1/b), & \text{if } \theta \in [0, C_1 b). \end{cases} \end{aligned}$$

PROOF. We begin by observing that, going back to (5.33),

$$
\begin{aligned}
|Q_2(z)| &= \left| \frac{1}{1-z} \left( \sum_{j=2}^{\infty} \frac{L(j)}{j^\alpha} \exp(-bj) + \sum_{j=2}^{\infty} \frac{L(j)}{j^{1+\alpha}} \exp(-bj) \sum_{n=0}^{j-1} z^n \right) \right| \\
&\leq \frac{2}{|1-z|} \sum_{j=2}^{\infty} \frac{L(j)}{j^\alpha} \exp(-bj) \leq 3\Gamma(1-\alpha) \frac{b^{\alpha-1} L(1/b)}{|1-z|},
\end{aligned}
$$
(5.53)

where, in the first equality, we assume that $z \neq 1$, as well as requiring that $|z| \leq 1$, and, in the second one, we assume $b$ to be sufficiently small. Going back to (5.44), we see that, for $|z| \leq 1$ and $b$ sufficiently small,

$$
\left| \frac{Q_2(z)}{Q_2(1)} \right| \leq \frac{6b}{(1-\alpha)|1-z|}.
$$
(5.54)

Note that this estimate becomes useless when $z$ is very close to 1.

On the other hand, with $z = \exp(-(B/2) + 2\pi\theta i)$ and $\widetilde{z} := \exp(2\pi\theta i)$, we write

$$
\begin{aligned}
(1-z)Q_1(z) &= \sum_{j=1}^{\infty} \frac{L(j)}{j^{1+\alpha}} \exp(-bj)(1-z^j) \\
&= \sum_{j=1}^{\infty} \frac{L(j)}{j^{1+\alpha}} (1-\widetilde{z}^j) - \sum_{j=1}^{\infty} \frac{L(j)}{j^{1+\alpha}} (1-\exp(-bj)) \\
&\quad + \sum_{j=1}^{\infty} \frac{L(j)}{j^{1+\alpha}} \widetilde{z}^j (1-\exp(-(b+(B/2))j)).
\end{aligned}
$$
(5.55)

Note that the absolute value of the last two terms is bounded above, respectively, by $b^\alpha L(1/b)$ and by $(b+(B/2))^\alpha L(1/(b+(B/2)))$ times a constant, which implies that their sum is bounded by $cb^\alpha L(1/b)$ (we assume, here, for example, that $B \leq b$) for some $c > 0$ and $b$ sufficiently small. Note, also, that the first term in the right-hand side of (5.55) is just $G(\theta)$; see (5.41). Since $|G(\cdot)| > 0$ is bounded away from 0 over compact subsets of $(0,1)$ (see Lemma 5.5) and given the asymptotic estimate (5.49), one directly sees that there exists a constant $C_1 > 0$ that guarantees that

$$
|1-z||Q_1(z)| \geq \tfrac{1}{2} |G(\theta)|
$$
(5.56)

if $1/2 \geq \theta/b \geq C_1$ and if $b$ is sufficiently small.

Therefore, using (5.37), we conclude that

$$
\left| \frac{Q_1(z)}{Q_1(1)} \right| \geq \frac{|G(\theta)| b^{1-\alpha}}{3|1-z|\Gamma(1-\alpha) L(1/b)}
$$
(5.57)

and putting (5.30), (5.57) together with (5.54), we obtain that, for $\theta \in [C_1 b, (1 - C_1 b)]$ and for $b$ sufficiently small,

$$(5.58) \qquad |\widehat{\Delta}_b(\exp(-(B/2) + 2\pi i \theta))| \leq \frac{10}{|G(\theta)|}.$$

Let us now turn our attention to $\theta \in (0, 1) \setminus [C_1 b, (1 - C_1 b)]$. In fact, it suffices to look at $\theta \in (0, C_1 b)$. For $F^{(2)}(\cdot)$, we will simply use the bound $|F^{(2)}(z)| \leq 1$, which holds if $|z| \leq 1$. So, we will just focus on finding a lower bound on $|F^{(1)}(\cdot)|$. It is technically convenient to separately consider the case of $\theta \in (0, \varepsilon_0 b]$ and $\theta \in [\varepsilon_0 b, C_1 b)$, where $\varepsilon_0 \in (0, C_1)$ is a small constant that we are going to choose below.

Let us start with the case $\theta \in (0, \varepsilon_0 b]$. In this case, note that $\Re \exp((-(B/2) + 2\pi i \theta)j) \geq 1/2$ if $\theta j \in [0, 1/8]$ and if $B$ is sufficiently small, so we can write

$$(5.59) \quad |Q_1(z)| \geq \frac{1}{2} \sum_{j \leq 1/(8\varepsilon_0 b)} \frac{L(j)}{j^\alpha} \exp(-bj) - \sum_{j > 1/(8\varepsilon_0 b)} \frac{L(j)}{j^\alpha} \exp(-bj)$$

and we now see that we can choose $\varepsilon_0$ so that

$$(5.60) \qquad |Q_1(z)| \geq \tfrac{1}{3} Q_1(1)$$

for $\theta \in (0, \varepsilon_0 b]$ and for every $b$ sufficiently small.

We are therefore left with the case of $\theta \in [\varepsilon_0 b, C_1 b)$. We are looking for a lower bound on the absolute value of $Q_1(\exp(-B/2) + 2\pi i \theta))$ and it is therefore sufficient to find a lower bound on the imaginary part of the same quantity. We use the elementary formula

$$(5.61) \qquad \Im\left(\frac{1 - z^j}{1 - z}\right) = \frac{\Im(1 - z^j)\Re(1 - \overline{z})}{|1 - z|^2} + \frac{\Re(1 - z^j)\Im(1 - \overline{z})}{|1 - z|^2}$$

and the fact that, in the regime we are in, as $b \searrow 0$, we have $1 - z \sim -2\pi i \theta$ so that $\Re(1 - \overline{z}) = o(b)$ and $\Im(1 - \overline{z}) \sim 2\pi \theta$. Let us also keep in mind that the ratio $\theta/b$ is bounded above and away from zero. The last considerations directly lead to the following two estimates when $z = \exp(-(B/2) + 2\pi i \theta)$: for every $\varepsilon > 0$, we have, for $b$ sufficiently small,

$$(5.62) \qquad \left|\frac{\Im(1 - z^j)\Re(1 - \overline{z})}{|1 - z|^2}\right| \leq \frac{\varepsilon}{2\pi\theta} |\sin(2\pi\theta j)| \leq \varepsilon j$$

and if $\theta j \leq 1/8$, one can find $c > 0$ such that

$$(5.63) \qquad \frac{\Re(1 - z^j)\Im(1 - \overline{z})}{|1 - z|^2} \geq c\theta j^2.$$

Therefore, for $b$ sufficiently small, considering that the second term in the right-hand side of (5.61) is nonnegative,

$$(5.64) \quad \begin{aligned} &|Q_1(\exp(-(B/2) + 2\pi i \theta))| \\ &\qquad \geq c\theta \sum_{j: \theta j \leq 1/8} L(j) j^{1-\alpha} \exp(-j\theta/\varepsilon_0) - \varepsilon \sum_j \frac{L(j)}{j^\alpha} \exp(-j\theta/C_1) \end{aligned}$$

and the right-hand side behaves, for $\theta \searrow 0$, like

$$(5.65) \quad L(1/\theta)\theta^{\alpha-1}\left(c\int_0^{1/8} t^{1-\alpha}\exp(-t/\varepsilon_0)\,\mathrm{d}t - \varepsilon\int_0^\infty t^{-\alpha}\exp(-t/C_1)\,\mathrm{d}t\right).$$

Since, if $\varepsilon$ is sufficiently small, the term between parentheses is positive, going back to (5.37), we see that there exists a positive constant $c'$ such that

$$(5.66) \quad \left|\frac{Q_1(\exp(-(B/2)+2\pi i\theta))}{Q_1(1)}\right| \ge c'$$

for every sufficiently small $b$. This concludes the proof of Lemma 5.6. □

5.1. *On the inter-arrival law of the intersection of renewals.* In this subsection, we study the asymptotic behavior of $\mathbb{K}_b(\cdot)$ itself. The case $b=0$ can be treated by using Banach space techniques (see [11] and references therein) and, while we could not find a precise reference in the literature to the specific estimate we needed (i.e., Proposition 5.7), we believe it is a classical result. When, instead, $b>0$, the outcome is somewhat surprising.

*The transient case.* We deal with the $b=0$ case and we are assuming (as usual) (2.11).

PROPOSITION 5.7. *We have*

$$(5.67) \qquad\qquad \mathbb{K}_0(n) \overset{n\to\infty}{\sim} c u_0(n)^2$$

*with $c=(\sum_{j=0}^\infty u_0(j)^2)^{-2} \in (0,1)$.*

PROOF. Let us set $P(n):=\mathbb{U}_0(n)/\widehat{\mathbb{U}}_0(1)$. Note that $\widehat{\mathbb{U}}_0(1)(=\sum_n u_0(n)^2) < \infty$ by (2.3) and (2.11). $P(\cdot)$ is therefore a probability distribution and we can write

$$(5.68) \qquad\qquad \widehat{\mathbb{K}}_0(z) = \phi(\widehat{P}(z))$$

with $\phi(\zeta) = 1 - (\widehat{\mathbb{U}}_0(1)\zeta)^{-1}$ and where $z$ is a complex number with $|z| \le 1$. In fact, 1 is the radius of convergence of the power series $\widehat{P}(\cdot)$. Notice that $\phi : \mathbb{C}\setminus\{0\} \longrightarrow \mathbb{C}$ is analytic and that there is no $z$ such that $\widehat{P}(z) = 0$ for $|z| \le 1$. This follows from (5.68) itself since $|\widehat{\mathbb{K}}_0(z)| \le \widehat{\mathbb{K}}_0(1) \le 1$ and a solution to $\widehat{P}(z) = 0$ in the unit ball would imply that $|\widehat{\mathbb{K}}_0(z)| = \infty$. Therefore, $\phi(\cdot)$ is analytic in a region (open connected set) containing the range of the power series $\widehat{P}(\cdot)$, that is, containing $\{\widehat{P}(z):|z|\le 1\}$. In this framework, one can apply Theorem 1 of [11] if some regularity properties on $P(\cdot)$ are verified. The regularity properties follow directly from the sufficient conditions in [11], page 259. The net outcome is that $\mathbb{K}_0(n) \overset{n\to\infty}{\sim} \phi'(\widehat{P}(1))P(n)$ and, since $\widehat{P}(1) = 1$, the result follows. □

*The positive recurrent case.* As mentioned before, in this regime, the results are somewhat unexpected. We will not strive for results in the most general setup, both because the results are mentioned merely as a side remark with respect to the main thrust of the paper and because complete answers are not obvious.

What we prove is essentially summed up by the following:

PROPOSITION 5.8. *Given $\alpha \in (0, 1/2)$ and the slowly varying function $L(\cdot)$, there exists $b_0 > 0$ such that, for $b \in (0, b_0)$, there exists $r \in (1, \exp(b))$ such that $(r-1)\widehat{\mathbb{D}}_b(r) = u_b^2(\infty)$. Moreover, if we call $r(b)$ the minimal value of $r$, then $r(b) > 1$ and $\mathbb{K}_b(n) \overset{n \to \infty}{\sim} C(n)r(b)^{-n}$ with $C(\cdot)$ such that $\limsup_n |C(n)| > 0$ and $|C(\cdot)|$ has at most polynomial growth.*

*If, in addition, $(z-1)\widehat{\mathbb{D}}_b(z) = u_b^2(\infty)$ has only one solution on $|z| = r(b)$ and this solution is a simple root, then there exists $C > 0$ such that*

$$(5.69) \qquad \mathbb{K}_b(n) \overset{n \to \infty}{\sim} Cr(b)^{-n}.$$

It will be clear from the argument below that once $r(b)$ is known, as well as any root of the equation in the statement on the circle of radius $r(b)$ centered at $0$, $C$ or, in the most general case, $C(\cdot)$ can be written explicitly.

Proposition 5.8 is somewhat surprising, particularly from a purely probabilistic viewpoint. By intersecting two independent renewals with inter-arrival laws such that $K_b(n) \exp(bn)$ is a nontrivial regularly varying function which vanishes at infinity, one may end up with inter-arrival laws that are purely exponential. As with Proposition 5.7, the proof involves complex analysis arguments and tail behaviors are linked to the different natures of the singularities that determine the radius of convergence of the $z$-transform of the sequence one studies.

PROOF OF PROPOSITION 5.8. Let us first recall that, for $b$ sufficiently small, we have $u_b(n) - u_b(\infty) \overset{n \to \infty}{\sim} (c(b)-1)^{-2}K_b(n)$ [19] [$c(b)$ is defined in (2.39)]. Note that this guarantees that the radius of convergence of $\widehat{\mathbb{D}}_b(z)$ is $\exp(b)$ and that $\widehat{\mathbb{D}}_b(z)$ converges at its radius of convergence. We have the formula

$$(5.70) \qquad \widehat{\mathbb{K}}_b(z) = 1 - \frac{1-z}{u_b^2(\infty) + (1-z)\widehat{\mathbb{D}}_b(z)}$$

for every $z$ in the centered ball of radius $\exp(b)$. Note that $\widehat{\mathbb{K}}_b(\cdot)$ is only meromorphic and we are actually going to prove that the denominator in the right-hand side does take the value zero. In order to see this, observe that, since $\widehat{\mathbb{K}}_b(\cdot)$ is a power series with nonnegative coefficients, if $\widehat{\mathbb{K}}_b(z)$ diverges, then $\widehat{\mathbb{K}}_b(|z|) = +\infty$. Let us concentrate on the real axis (there may also be complex poles with the same absolute value). We claim that there is at least one singularity at $z = r \in (1, \exp(b))$ for $b$ sufficiently small. To see this, note that the denominator in the right-hand side of (5.70) is analytic in the centered ball of radius $\exp(b)$ and it takes the values $u_b^2(\infty)$ at $z = 1$. At $z = \exp(b)$, it instead takes the value $u_b^2(\infty) - (\exp(b) - 1)\widehat{\mathbb{D}}_b(\exp(b))$. From (5.24), we have

$$(5.71) \qquad \widehat{\mathbb{D}}_b(\exp(b)) \geq 2u_b(\infty)\widehat{\Delta}_b(\exp(b)).$$

To evaluate $\widehat{\Delta}_b(\exp(b))$, we observe that

$$(5.72) \qquad \widehat{\Delta}_b(\exp(b)) = \frac{1}{1 - c(b)} - \frac{u_b(\infty)}{1 - \exp(b)} \overset{b \searrow 0}{\sim} (1 - \alpha)\frac{u_b(\infty)}{b}$$

so that the denominator in the right-hand side of (5.70) is negative for $z = \exp(b)$ and $b$ sufficiently small; in fact, it is asymptotically equivalent to $(1 - 2(1 - \alpha))u_b(\infty)^2$ (use the estimates in Remark 5.1). This implies that there exists some $z \in (1, \exp(b))$ for which the denominator in the right-hand side of (5.70) takes the value zero: $z$ is called $r$ in the statement. Note, also, that any solution $r$ is necessarily bounded away from 1.

The rest of the proof is just based on standard expansions at the pole singularities of $\widehat{\mathbb{K}}_b(\cdot)$ that are closest to the origin; see, for example, [19], Section 4.2. $\square$

**Note added in proof.** After the completion of this work, a number of results on the relevant disorder regime of pinning models have been proven by K. Alexander, B. Derrida, H. Lacoin, N. Zygouras and by the authors of this work (some of these results and an updated bibliography may be found in [20]).

UNIVERSITÉ PARIS DIDEROT (PARIS 7)
  AND LABORATOIRE DE PROBABILITÉS
  ET MODÈLES ALÉATOIRES (CNRS)
U.F.R. MATHÉMATIQUES
CASE 7012 (SITE CHEVALERET)
75205 PARIS CEDEX 13
FRANCE
E-MAIL: giacomin@math.jussieu.fr

CNRS AND LABORATOIRE DE PHYSIQUE
ENS LYON
46 ALLÉE D'ITALIE
69364 LYON
FRANCE
E-MAIL: fabio-lucio.toninelli@ens-lyon.fr